\documentclass[leqno,draft,11pt]{amsart}
\input{amssym.def}

\newcommand{\RR}{{\Bbb R}}
\newcommand{\NN}{{\Bbb N}}

\newcommand{\fin}{\hfill$\Box$}

\newtheorem{teo}[equation]{Theorem}

\newtheorem{lema}[equation]{Lemma}
\newtheorem{propo}[equation]{Proposition}

\theoremstyle{definition} \theoremstyle{remark}

\numberwithin{equation}{section}

\begin{document}

\title[]
{Higher order Riesz transforms for Laguerre expansions}

 \keywords{}
\begin{abstract}
In this paper we investigate $L^p$-boundedness properties for the
higher order Riesz transforms associated with Laguerre operators.
Also we prove that the $k$-th Riesz transform is a principal value
singular integral operator (modulus a constant times of the function
when $k$ is even). To establish our results we exploit a new
identity connecting Riesz transforms in the Hermite and Laguerre
settings.
\end{abstract}

\author[J. Betancor]{Jorge J. Betancor}
\address{Departamento de An\'{a}lisis Matem\'{a}tico\\
Universidad de la Laguna\\
Campus de Anchieta, Avda. Astrof\'{\i}sico Francisco S\'{a}nchez, s/n\\
38271 La Laguna (Sta. Cruz de Tenerife), Spain}
\email{jbetanco@ull.es; jcfarina@ull.es; lrguez@ull.es;
asgarcia@ull.es}

\author[J.C. Fari\~{n}a]{Juan C. Fari\~{n}a}

\author[L. Rodr\'{\i}guez-Mesa]{Lourdes Rodr\'{\i}guez-Mesa}

\author[A. Sanabria-Garc\'{\i}a]{Alejandro Sanabria-Garc\'{\i}a}

\thanks{This paper is partially supported by MTM2007/65009. Third and fourth authors are
also partially supported by grant PI042004/067}
\date{\today}

\maketitle

\section{Introduction}

The aim of this paper is to investigate higher order Riesz
transforms associated with Laguerre function expansions. To
achieve our goal we use a procedure that will be described below
and that was developed for the first time by the authors and
Torrea in \cite{BFRST}. Our results complete and improve in some
senses the ones obtained by Graczyk, Loeb, L\'opez, Nowak and
Urbina \cite{GLLNU} about higher order Riesz transforms for
Laguerre expansions.

For every $\alpha >-1$ we consider the Laguerre differential
operator
$$
L_\alpha
=\frac{1}{2}\left(-\frac{d^2}{dx^2}+x^2+\frac{1}{x^2}\left(\alpha
^2-\frac{1}{4}\right)\right),\quad x\in (0,\infty ).
$$

This operator can be factorized as follows
\begin{equation}\label{factor}
L_\alpha =\frac{1}{2}\mathfrak{D}_\alpha ^*\mathfrak{D}_\alpha
+\alpha +1,
\end{equation}
where $\mathfrak{D}_\alpha f=\left(\displaystyle -\frac{\alpha
+1/2}{x}+x+\frac{d}{dx}\right)f=x^{\alpha +\frac{1}{2}}\displaystyle
\frac{d}{dx}(x^{-\alpha -\frac{1}{2}}f)+xf$, and
$\mathfrak{D}_\alpha ^*$ denotes the formal adjoint of
$\mathfrak{D}_\alpha$ in $L^2((0,\infty ),dx)$.

The heat semigroup $\{W_t^\alpha\}_{t>0}$ generated by the Laguerre
operator $-L_\alpha$ admits the integral representation
\begin{equation}
W_t^\alpha (f)(x)=\int_0^\infty W_t^\alpha (x,y)f(y)dy,\quad x\in
(0,\infty ),\; f\in L^2((0,\infty ),dx)\;, \label{integral}
\end{equation}
 where the heat kernel $W_t^\alpha(x,y)$, $t,x,y\in
(0,\infty)$, is defined in (\ref{ker1}). The factorization
(\ref{factor}) for $L_\alpha$ suggests to define (formally), for
every $k\in \mathbb{N}$, the $k$-th Riesz transform $R_\alpha^{(k)}$
associated with $L_\alpha$ by
$$
R_\alpha^{(k)}=\mathfrak{D}_\alpha^k L_\alpha ^{-\frac{k}{2}}\;.
$$
Here $L_\alpha^{-\beta}$, $\beta>0$, denotes the $-\beta$ power of
the operator $L_\alpha$ (see (\ref{Pow1})).

In the main result of this paper (Theorem \ref{main}) we prove that
the space $C_c^\infty(0,\infty)$ of $C^\infty$ functions on
$(0,\infty)$ having compact support on $(0,\infty)$ is contained in
the domain of $R_\alpha^{(k)}$ and that $R_\alpha^{(k)}$ on
$C_c^\infty(0,\infty)$ is a principal value integral operator.
Moreover we establish that $R_\alpha^{(k)}$ can be extended as a
principal value integral and bounded operator on certain weighted
$L^p$ -spaces.

\begin{teo}\label{main}
Let $\alpha >-1$ and $k\in \NN$. For every  $\phi \in C_c^\infty
(0,\infty )$ it has that
$$
R_\alpha ^{(k)}\phi (x)=w_k\phi (x)+\lim_{\varepsilon \rightarrow
0^+}\int_{0,|x-y|>\varepsilon }^\infty R_\alpha ^{(k)}(x,y)\phi
(y)dy,\quad x\in (0,\infty ),
$$
where
$$
R_\alpha ^{(k)}(x,y)=\frac{1}{\Gamma
\Big(\frac{k}{2}\Big)}\int_0^\infty t^{\frac{k}{2}-1}
\mathfrak{D}_\alpha ^kW_t^\alpha (x,y)dt,\quad x,y\in (0,\infty ),
$$
and $w_k=0$, when $k$ is odd and $w_k=-2^{\frac{k}{2}}$, when $k$ is
even.

The operator $R_\alpha ^{(k)}$ can be extended, defining it by
$$
R_\alpha ^{(k)}f(x)=w_kf(x)+\lim_{\varepsilon \rightarrow
0^+}\int_{0,|x-y|>\varepsilon }^\infty R_\alpha
^{(k)}(x,y)f(y)dy,\quad \mbox{ a.e. } x\in (0,\infty ),
$$
as a bounded operator from $L^p((0,\infty ), x^\delta dx)$ into
itself, for $1<p<\infty$ and

(a) $-\Big(\alpha +\frac{3}{2}\Big)p-1<\delta < \Big(\alpha
+\frac{3}{2}\Big)p-1$, when $k$ is odd;

(b) $-\Big(\alpha +\frac{1}{2}\Big)p-1<\delta < \Big(\alpha
+\frac{3}{2}\Big)p-1$, when $k$ is even;

\noindent and as a bounded operator from $L^1((0,\infty ), x^\delta
dx)$ into $L^{1,\infty }((0,\infty ), x^\delta dx)$ when

(c) $-\alpha -\frac{5}{2}\leq \delta \leq \alpha +\frac{1}{2}$, when
$k$ is odd;

(d) $ -\alpha -\frac{3}{2}\leq \delta \leq \alpha +\frac{1}{2}$, for
$\alpha \not=-\frac{1}{2}$, and $-1<\delta \leq 1$, for $\alpha
=-\frac{1}{2}$, when $k$ is even.
\end{teo}

Also we get the corresponding property in the Hermite context (see
Proposition \ref{Hermite}) which completes in the one dimensional
case the results in \cite{StTo2} about the higher order Riesz
transform associated with the Hermite operator.

For every $n\in \NN$, we have that $L_\alpha \varphi _n^\alpha
=(2n+\alpha +1)\varphi _n^\alpha$, where,
$$
\varphi _n^\alpha (x)=\left(\frac{2\Gamma (n+1)}{\Gamma (n+\alpha
+1)}\right)^{\frac{1}{2}}e^{-\frac{x^2}{2}}x^{\alpha
+\frac{1}{2}}L_n^\alpha (x^2),\quad x\in (0,\infty ),
$$
and $L_n^\alpha$ denotes the $n$-th Laguerre polynomial of type
$\alpha$ (\cite[p. 100]{Szeg} and \cite[p. 7]{Than}). For every
$n\in \NN$, $\varphi _n^\alpha $ is usually called the $n$-th
Laguerre function of type $\alpha$. The system $\{\varphi
_n^\alpha \}_{n\in \NN}$ is an orthonormal basis for
$L^2((0,\infty ),dx)$.

The heat semigroup $\{W_t^\alpha \}_{t>0}$ generated by the operator
$-L_\alpha$ is defined by
$$
W_t^\alpha (f)=\sum_{n=0}^\infty e^{-t(2n+\alpha +1)}c_n^\alpha
(f)\varphi _n^\alpha ,\quad f\in L^2((0,\infty ),dx),
$$
where $c_n^\alpha (f)=\displaystyle \int_0^\infty \varphi _n^\alpha
(x)f(x)dx$, $n\in \NN$.

For every $t>0$ the operator $W_t^\alpha$ admits the integral
representation (\ref{integral}) where, by the Mehler formula
\cite[p. 8]{Than}, for every $t,x,y\in (0,\infty )$,
\begin{eqnarray}
W_t^\alpha (x,y)&=&\displaystyle \sum_{n=0}^\infty e^{-t(2n+\alpha
+1)}\varphi _n^\alpha (x)\varphi _n^\alpha (y)\vspace{3mm}\nonumber\\
&=&\displaystyle
\left(\frac{2e^{-t}}{1-e^{-2t}}\right)^{\frac{1}{2}}
\left(\frac{2xye^{-t}}{1-e^{-2t}}\right)^{\frac{1}{2}}I_\alpha
\left(\frac{2xye^{-t}}{1-e^{-2t}}\right)e^{-\frac{1}{2}(x^2+y^2)
\frac{1+e^{-2t}}{1-e^{-2t}}}. \label{ker1}
\end{eqnarray}
Here $I_\alpha$ represents the modified Bessel function of the first
kind and order $\alpha$.

According to the ideas presented by Stein \cite{Stei} the
fundamental operators of the harmonic analysis (fractional
integrals, Riesz transforms, g-functions,...) can be considered in
the context of the Laguerre operator $L_\alpha$. It is convenient
to mention that this way to describe harmonic operators in the
Laguerre context was initiated by Muckenhoupt (\cite{Muck1} and
\cite{Muck3}).

If $\beta >0$ the negative power $L_\alpha ^{-\beta}$ of $L_\alpha$
is defined by
$$
L_\alpha^{-\beta} f = \sum_{n=0}^\infty \frac{c_n^\alpha
(f)}{(2n+\alpha+1)^\beta} \varphi_n^\alpha, \quad f \in
L^2((0,\infty),dx).
$$

It is not hard to see that $L_\alpha ^{-\beta}$ can be expressed,
for every $f\in L^2((0,\infty ),dx)$, by means of the following
integral
\begin{equation}
L_\alpha^{-\beta} f (x)= \frac{1}{\Gamma (\beta )}\int_0^\infty
t^{\beta -1}W_t^\alpha (f)(x)dt,\quad \mbox{ a.e. }x\in (0,\infty ).
\label{Pow1}
\end{equation}

$L_\alpha ^{-\beta}$ is also called $\beta$-th fractional integral
associated with $L_\alpha$. This kind of fractional integrals has
been investigated by several authors (\cite{ChiPola2},
\cite{GaStTr}, \cite{GLLNU}, \cite{Kanj} and \cite{Stem}).

First order Riesz transforms in the $L_\alpha$-setting were studied
in \cite{NoSt2} for $\displaystyle \alpha \geq -\frac{1}{2}$ and in
\cite{AMST} for $\alpha >-1$. Also, the procedure developed in
\cite{HSTV} can be used to investigate strong, weak and restricted
weak type with respect to the measure $x^\delta dx$ on $(0,\infty )$
for the Riesz transforms $R_\alpha^{(1)}$.

As it was mentioned, in this paper we establish boundedness
properties for $R_\alpha ^{(k)}$ in $L^p((0,\infty), x^\delta dx)$.
In the Laguerre-polynomial context, Graczyk, Loeb, L\'opez, Nowak
and Urbina \cite{GLLNU} investigated the corresponding higher order
Riesz transform $\mathcal{R}_\alpha ^{(k)}$, $k\in \NN$. They used
methods that impose substancial restrictions to the admisible values
of the index $\alpha$. Their proofs are based in a connection
between $n$-dimensional Hermite functions and Laguerre functions of
order $\alpha =\displaystyle \frac{n}{2}-1$. This fact forces to
consider only half-integer values for $\alpha$ (see \cite[Section
4]{GLLNU}). This connection between $n$-dimensional Hermite and
Laguerre functions was exploited earlier by Guti\'errez, Incognito
and Torrea \cite{GuInTo} and Harboure, Torrea and Viviani
\cite{HaToVi}, amongst others. As it is emphatized in \cite{GLLNU},
the extension of the results of $L^p$-boundedness for
$\mathcal{R}_\alpha ^{(k)}$ to all values of $\alpha
>-1$ is not an easy problem and it requires a more subtle
approach than the one followed in \cite{GLLNU}.

Our procedure here is completely different from the one used in
\cite{GLLNU}. In a first step we split the operators $R_\alpha
^{(k)}$ into two parts, namely: a local operator and a global one.
These operators are integral operators defined by kernels
supported close to and far from the diagonal, respectively. The
global operator is upper bounded by Hardy type operators. The
novelty of our method is the way followed to study the local part.
We establish a pointwise identity connecting the kernel of
$R_\alpha ^{(k)}$ with the one corresponding to the $k$-th Riesz
transform associated with the Hermite operator in one dimension,
for every $\alpha >-1$ (see Proposition \ref{kernels}). By using
this identity we transfer boundedness and  convergence results
from $k$-th Riesz transform for Hermite operator in one dimension
to $k$-th Riesz transform in the $L_\alpha$-setting.

In the literature (see, for instance \cite{ChiPola} and \cite{Stem})
we can find other systems of Laguerre functions different from
$\{\varphi _n^\alpha \}_{n\in \NN}$. In particular, from the
Laguerre polynomials $\{L_n^\alpha \}_{n\in \NN}$ we can derive also
the systems $\{\mathcal{L}_n^\alpha \}_{n\in \NN}$ and $\{l_n^\alpha
\}_{n\in \NN}$, where, for every $n\in \NN$,
$$
\mathcal{L}_n^\alpha
(x)=\left(\frac{\Gamma (n+1)}{\Gamma (n+\alpha
+1)}\right)^{\frac{1}{2}}e^{-\frac{x}{2}}x^{\frac{\alpha
}{2}}L_n^\alpha (x),\quad x\in (0,\infty ),
$$
and
$$
l_n^\alpha (x)=\left(\frac{\Gamma (n+1)}{\Gamma (n+\alpha
+1)}\right)^{\frac{1}{2}}e^{-\frac{x}{2}}L_n^\alpha (x),\quad x\in
(0,\infty )\;.
$$
$\{\mathcal{L}_n^\alpha \}_{n\in \NN}$ is an orthonormal basis in
$L^2((0,\infty ),dx)$ and $\{l_n^\alpha \}_{n\in \NN}$ is an
orthonormal basis in $L^2((0,\infty ), x^\alpha dx)$.

As it is shown in \cite{AMST}, harmonic analysis operators
associated with $\{\mathcal{L}_n^\alpha \}_{n\in \NN}$ and
$\{l_n^\alpha \}_{n\in \NN}$ is closely connected with the
corresponding operators related to the family $\{\varphi _n^\alpha
\}_{n\in \NN}$. The connection is given by a multiplication
operator defined by $M_\beta f=x^\beta f$, for certain $\beta \in
\RR$. From the strong type results for $R_\alpha ^{(k)}$
established in Theorem \ref{main}, the corresponding results for
the $k$-th Riesz transform in the $\{\mathcal{L}_n^\alpha \}_{n\in
\NN}$ and $\{l_n^\alpha \}_{n\in \NN}$ settings can be deduced.
Moreover the weak type results for the $k$-th Riesz transform
associated with $\{\mathcal{L}_n^\alpha \}_{n\in \NN}$ and
$\{l_n^\alpha \}_{n\in \NN}$ can be obtained by proceeding as in
the $\{\varphi _n^\alpha \}_{n\in \NN}$ case of Theorem
\ref{main}.

The organization of the paper is the following. Section 2 contains
some basic facts needed in the sequel. Section 3 is devoted to prove
the main result of this paper (Theorem \ref{main}) where we
establish that the higher order Riesz transforms are principal value
singular integral operators (modulus a constant times of the
function, when $k$ is even) and we show $L^p((0,\infty ),x^\delta
dx)$-boundedness properties for them.

Throughout this paper $C_c^\infty (I)$ denotes the space of
functions in $C^\infty (I)$ having compact support on $I$. By $C$
and $c$ we always represent positive constants that can change
from one line to the other one, and $E[r]$, $r\in \RR$, stands for
the integer part of $r$.

The authors would like to express their gratitude to Prof. Jos\'e
Luis Torrea for its valuable and helpful comments which have
certainly led to the improvement of the paper.

\section{Preliminaries}
In this section we recall some definitions and properties that will
be useful in the sequel. By $H$ we denote the Hermite differential
operator
\begin{equation}\label{Hermiteop}
H=\frac{1}{2}\left(-\frac{d^2}{dx^2}+x^2\right)=-\frac{1}{4}\left[\left(\frac{d}{dx}+x\right)
\left(\frac{d}{dx}-x\right)+\left(\frac{d}{dx}-x\right)\left(\frac{d}{dx}+x\right)\right]\;.
\end{equation}
Note that $\frac{d}{dx}+x$ and $-\frac{d}{dx}+x$ are formal
adjoint operators in $L ^2(\RR ,dx)$. Moreover, if $n\in \NN$,
$H_n$ represents the $n$-th Hermite polynomial (\cite[p.
104]{Szeg}) and $h_n$ is the Hermite function given by $
h_n(x)=(\sqrt{\pi }2^nn!)^{-\frac{1}{2}}e^{-\frac{x^2}{2}}H_n(x)$,
$ x\in \RR$, then it has that
$$
Hh_n=\left(n+\frac{1}{2}\right)h_n,\quad n\in \NN \;.
$$
Moreover, the system $\{h_n\}_{n\in \NN}$ is an orthonormal basis
in $L^2(\RR ,dx)$.

The investigations of harmonic analysis in the Hermite setting
were begun by Muckenhoupt (\cite{Muck1}). This author considered
Hermite polynomial expansions instead of Hermite function
expansions. In the last decades several authors have studied
harmonic analysis operators in the Hermite (polynomial or
function) context (see, for instance, \cite{FaGuSc}, \cite{FoSc},
\cite{GaMaSjTo}, \cite{GuSeTo}, \cite{HaRoSeTo}, \cite{PeSo},
\cite{Sjog}, \cite{StTo1}, \cite{StTo2} and \cite{Urbi}).

The heat semigroup $\{W_t\}_{t>0}$ generated by the operator $-H$ is
defined by
$$
W_t(f)=\sum_{n=0}^\infty e^{-(n+\frac{1}{2})t}c_n(f)h_n,\quad f\in
L^2(\RR, dx)\mbox{ and }t>0,
$$
where $c_n(f)=\displaystyle \int_{-\infty }^{+\infty}
h_n(x)f(x)dx$, $n\in \NN$ and $f\in L^2(\RR ,dx)$.

For every $t>0$, the operator $W_t$ can be described by the
integral
$$
W_t(f)(x)=\int_{-\infty}^{+\infty} W_t(x,y)f(y)dy,\quad x\in \RR
\mbox{ and }f\in L^2(\RR ,dx),
$$
where, according to the Mehler formula \cite[p. 2]{Than}, we have
that, for each $x,y\in \RR$ and $t>0$,
$$
W_t(x,y)=\sum_{n=0}^\infty
e^{-(n+\frac{1}{2})t}h_n(x)h_n(y)=\frac{1}{\sqrt{\pi
}}\left(\frac{e^{-t}}{1-e^{-2t}}\right)^{\frac{1}{2}}e^{-\frac{1}{2}(x^2+y^2)
\frac{1+e^{-2t}}{1-e^{-2t}}+\frac{2xye^{-t}}{1-e^{-2t}}}.
$$

The negative power $H^{-\beta}$, $\beta >0$, of $H$ is given by
$$
H^{-\beta }f=\sum_{n=0}^\infty
\frac{c_n(f)}{\Big(n+\frac{1}{2}\Big)^\beta}h_n,\quad f\in L^2(\RR
,dx).
$$
It can be seen that, for every $f\in L^2(\RR ,dx)$,
$$
H^{-\beta }f(x)=\frac{1}{\Gamma (\beta )}\int_{-\infty}^{+\infty}
t^{\beta -1}W_t(f)(x)dt,\quad \mbox{ a.e. }x\in \RR \;,
$$
and that, when $\phi \in C_c^\infty (\RR )$,
$$
H^{-\beta }\phi (x)=\int_{-\infty}^{+\infty} K_{2\beta} (x,y)\phi
(y)dy,\quad x\in \RR \;.
$$
Here, for every $\gamma >0$,
$$
K_\gamma (x,y)=\frac{1}{\Gamma (\frac{\gamma }{2})}\int_0^\infty
t^{\frac{\gamma }{2} -1}W_t(x,y)dt,\quad x,y\in \RR \;.
$$

The factorization in (\ref{Hermiteop}) suggests to define the Riesz
transform $R$ associated with $H$ by
$$
Rf=\left(\frac{d}{dx}+x\right)H^{-\frac{1}{2}}f=\sum_{n=1}^\infty
\left(\frac{2n}{n+\frac{1}{2}}\right)^{\frac{1}{2}}c_n(f)h_{n-1},\quad
f\in L^2(\RR ,dx)\;.
$$
The operator $R$ admits the integral representation
$$
Rf(x)=\int_{-\infty}^{+\infty} R(x,y)f(y)dy,\quad x\in \RR \setminus
\mbox{supp } f\mbox{ and }f\in L^2(\RR ,dx),
$$
where
$$
R(x,y)=\frac{1}{\sqrt{\pi }}\int_0^\infty
t^{-\frac{1}{2}}\Big(\frac{d}{dx}+x\Big)W_t(x,y)dt,\quad x,y\in \RR
\;.
$$
$L^p$-boundedness properties of the Riesz transform $R$ (even in
the $n$-dimensional case) were established in \cite{StTo1}.

In \cite{StTo2} higher order Riesz transforms in the Hermite
function setting on $\RR ^n$ were investigated. For every $k\in \NN$
the $k$-th Riesz transforms on $\RR$, $R^{(k)}$, is defined by
$$
R^{(k)}f=\Big(\frac{d}{dx}+x\Big)^kH^{-\frac{k}{2}}f=\sum_{n=k}^\infty
\frac{2^{\frac{k}{2}}(n(n-1)...(n-k+1))^{\frac{1}{2}}}{(n+\frac{1}{2})
^{\frac{k}{2}}}c_n(f)h_{n-k},\quad f\in C_c^\infty (\RR ).
$$
$L^p$-boundedness properties of the Riesz transform $R^{(k)}$ were
established in \cite[Theorem 2.3]{StTo2} by invoking
Calder\'on-Zygmund singular integral theory. For every $k\in \NN$,
$R^{(k)}$ can be extended to $L^p(\RR ,dx)$ as a bounded operator
from $L^p(\RR ,dx)$ into itself, when $1<p<\infty$, and from
$L^1(\RR ,dx)$ into $L^{1,\infty }(\RR ,dx)$. It is remarkable to
note that the $L^p$-mapping properties for the higher order Riesz
transform in the Hermite polynomial setting are essentially
different to the corresponding ones in the Hermite function
context (\cite{FoSc} and \cite{GaMaSjTo}).

We shall not use the Calder\'on-Zygmund singular integral theory
to investigate the higher order Riesz transforms associated with
the Laguerre operator. As it was mentioned we shall exploit a
connection between higher order Riesz transforms in the Hermite
and Laguerre settings.

The following new property that will be established in Section 3 is
needed in the proof of Theorem \ref{main}. It states that the higher
order Riesz transform for the Hermite operator is actually a
principal value integral operator.

\begin{propo}\label{Hermite}
Let $k\in \NN$. Then, for every $f\in C_c^\infty(\mathbb{R})$,
$1\leq p<\infty $,
\begin{equation}
R^{(k)}f(x)=w_kf(x)+\lim_{\varepsilon \rightarrow
0^+}\int_{|x-y|>\varepsilon }R^{(k)}(x,y)f(y)dy,\quad \mbox{ a.e.
}x\in \RR , \label{vp2}
\end{equation}
where
$$
R^{(k)}(x,y)=\frac{1}{\Gamma(\frac{k}{2})}\int_0^\infty
t^{\frac{k}{2}-1}\left(\frac{d}{dx}+x\right)^kW_t(x,y)dt,\quad
x,y\in \RR ,
$$
and $w_k=0$, when $k$ is odd, and $w_k=-2^{\frac{k}{2}}$, when $k$
is even.
\end{propo}

Since $R^{(k)}(x,y)$, $x,y\in \mathbb{R}$, is a Calder\'on-Zygmund
kernel (\cite{StTo2}) by using standard density arguments from
Proposition \ref{Hermite} we deduce that the operator $R^{(k)}$ can
be extended by (\ref{vp2}) to $L^p(\mathbb{R},dx)$, $1\le p<\infty$,
as a bounded operator from $L^p(\mathbb{R},dx)$ into itself,
$1<p<\infty$, and from $L^1(\mathbb{R},dx)$ into
$L^{1,\infty}(\mathbb{R},dx)$.

As it was indicated the modified Bessel function $I_\alpha$ of the
first kind and order $\alpha$ appears in the kernel of the heat
semigroup generated by the Laguerre operator $-L_\alpha$. The
following properties of the function $I_\alpha$ will be repeteadly
used in the sequel (see \cite{Lebe} and \cite{Wats}):

\noindent (P1) $I_\alpha (z)\sim z^\alpha$, $z\rightarrow 0$.

\noindent (P2) $\displaystyle \sqrt{z}I_\alpha
(z)=\frac{e^z}{\sqrt{2\pi }}\left(\sum_{r=0}^n(-1)^r[\alpha ,
r](2z)^{-r}+O(z^{-n-1})\right)$, $n=0,1,2,...$, where $[\alpha
,0]=1$ and
$$
[\alpha,r] = \frac{(4\alpha^2-1) (4\alpha^2-3^2) \cdots
(4\alpha^2-(2r-1)^2)}{2^{2r} \Gamma(r+1)},\quad r=1,2,...\;.
$$

\noindent (P3) $\displaystyle \frac{d}{dz}(z^{-\alpha }I_\alpha
(z))=z^{-\alpha }I_{\alpha +1}(z)$, $z\in (0,\infty )$.

On the other hand, in our study for the global part of the
operators we consider the Hardy type operators defined by
$$
H_0^\eta (f)(x)=x^{-\eta -1}\int_0^xy^\eta f(y)dy,\quad x\in
(0,\infty ),
$$
and
$$
H_\infty ^\eta (f)(x)=x^\eta \int_x ^{\infty} y^{-\eta
-1}f(y)dy,\quad x\in (0,\infty ),
$$
where $\eta >-1$. $L^p$-boundedness properties of the operators
$H_0^\eta$ and $H_\infty ^\eta$ were established by Muckenhoupt
\cite{Muck4} and Andersen and Muckenhoupt \cite{AnMu}. In
particular, mappings properties for $H_0^\eta$ and $H_\infty
^\eta$ on $L^p((0,\infty ),x^\delta dx)$ can be encountered in
\cite[Lemmas 3.1 and 3.2]{ChiPola}.

The following formula established in \cite[Lemma 4.3,
(4.6)]{GLLNU} will be frequently used in the sequel. For every
$N\in\mathbb{N}$, and a sufficiently smooth function $g: (0,\infty
)\longrightarrow \RR$, it has that
\begin{equation}\label{2.3}
\frac{d^N}{dx^N} [g(x^2)] = \sum_{l=0}^{E[\frac{N}{2}]} E_{N,l}
x^{N-2l} \left(\frac{d^{N-l}}{dx^{N-l}} g \right)(x^2)
\end{equation}
where
$$
E_{N,l} = 2^{N-2l} \frac{N!}{l!(N-2l)!}, \quad 0\le l \le
E\Big[\frac{N}{2}\Big].
$$

We finish this section establishing the following technical lemma
that is needed in the proof of Proposition \ref{kernels}.

\begin{lema}\label{technical}
Let $\alpha >-1$ and $j\in \NN$, $j\geq 1$. For every $m=0,1,
\dots, E[\frac{j}{2}]$, we have
\begin{equation}\label{N1}
\sum_{n=0}^m \sum_{l=2n}^j (-1)^{l+n}{j \choose l}
\frac{E_{l,n}}{2^{l-2n}} [\alpha +l-n, m-n] =0.
\end{equation}
\end{lema}

\noindent \underline {Proof}. For every $s=0,...,j$ we denote by
$A_{j,s}$ the values
$$
A_{j,s}=\sum_{l=0}^j (-1)^l {j \choose l} l^s,
$$
where we take the convention $0^0=1$.

In \cite[(43)]{BFMR} it was established that, for every
$j\in\mathbb{N}$, $j\ge 1$,
\begin{equation}\label{N2}
A_{j,s}=0, \quad s=0,1,\cdots ,j-1.
\end{equation}
On the other hand, since ${m\choose n } = \frac{m}{n} {{m-1}
\choose {n-1}}$, for $m\ge n \ge 1$, by using (\ref{N2}) we obtain
that $A_{j,j}=-jA_{j-1,j-1}$, $j\in \NN$, $j\geq 1$ and so
$A_{j,j}=(-1)^jj!$, $j\in \NN$.

Fix $j\in \NN$, $j\geq 1$ and consider $m=0,1,...,E[\frac{j}{2}]$
and $n=0,1,...,m$. We can write
\begin{eqnarray*}
\sum_{l=2n}^j (-1)^l{j \choose l}\frac{E_{l,n}}{2^{l-2n}}
[\alpha+l-n,m-n] &=& \sum_{s=0}^{j-2n} (-1)^s{{j}\choose {s+2n}}
\frac{E_{s+2n,n}}{2^s} [\alpha+s+n,m-n]\\
&=& \frac{j!}{n!} \sum_{s=0}^{j-2n} \frac{(-1)^s}{s!(j-2n-s)!}
[\alpha+s+n,m-n]\\
&=& \frac{j!}{n!(j-2n)!} \sum_{s=0}^{j-2n}(-1)^s{{j -2n} \choose s}
[\alpha+s+n,m-n]\;.
\end{eqnarray*}
We observe that $[\alpha+s+n,m-n]$ is a polynomial in $s$ which
has degree $2(m-n)$. Besides, if $j$ is odd, $2(m-n)\leq j-2n-1$.
Hence (\ref{N2}) allows us to conclude (\ref{N1}) in this case.

Assume now that $j$ is even. Then $2(m-n)\leq j-2n$ and (\ref{N2})
leads to
$$
\displaystyle \sum_{n=0}^m \sum_{l=2n}^j(-1)^{l+n} {j\choose l}
\frac{E_{l,n}}{2^{l-2n}}[\alpha +l-n,m-n]=0,
$$
when $m=0,1,...,\frac{j}{2}-1$.

For $m=\frac{j}{2}$ and again by (\ref{N2}) we can write
$$
\begin{array}{l}
\displaystyle \sum_{n=0}^m \sum_{l=2n}^j(-1)^{l+n} {j\choose l}
\frac{E_{l,n}}{2^{l-2n}}[\alpha +l-n,m-n]
=\sum_{n=0}^{\frac{j}{2}}\frac{(-1)^nj!}{n!(j-2n)!}\sum_{s=0}^{j-2n}
(-1)^s{{j-2n}\choose s}
\frac{s^{j-2n}}{(\frac{j}{2}-n)!}\\
\\
\;\;\displaystyle = j!\sum_{n=0}^{\frac{j}{2}}
\frac{(-1)^n}{n!(j-2n)!(\frac{j}{2} -n)!}A_{j-2n,j-2n}
=\frac{j!}{(\frac{j}{2})!}\sum_{n=0}^{\frac{j}{2}}
(-1)^n{{\frac{j}{2}} \choose n}=0\;.
\end{array}
$$

Thus (\ref{N1}) is established. \fin

\section{Higher order Riesz transforms associated with Laguerre
expansions}

In this section we prove our main result (Theorem \ref{main})
concerning to  higher order  Riesz transforms associated with the
sequence $\{\varphi_n^\alpha\}_{n\in \NN}$ of Laguerre functions.
As it was mentioned, our procedure is based on certain connection
between higher order Riesz transforms in the Laguerre and Hermite
settings.

We start proving the following results about the differentiability
of $H^{-\frac{k}{2}}$, $k\in \NN$. For each $k\in \NN$ and
$l=0,1,...,k$, let us consider
$$
R^{(k,l)}(x,y)=\frac{1}{\Gamma (\frac{k}{2})}\int_0^\infty
t^{\frac{k}{2}-1}\Big(\frac{d}{dx}+x\Big)^lW_t(x,y)dt,\quad x,y\in
\RR .
$$
(Note that $R^{(k,k)}(x,y)=R^{(k)}(x,y)$, $x,y\in \RR$).

\begin{propo}\label{derivH}
Let $\phi \in C_c^\infty (\RR )$ and $k\in \NN$. We have that
$$
\Big(\frac{d}{dx}+x\Big)^lH^{-\frac{k}{2}}\phi (x)=\int_{-\infty
}^{+\infty} R^{(k,l)}(x,y)\phi (y)dy,
$$
for every $x\in \RR$, when $l=0,...,k-1$ and for every $x\in \RR
\setminus {\rm supp}\;\phi $, when $l=k$.
\end{propo}

\noindent \underline{Proof}. As it was mentioned, we can write
$$
H^{-\frac{k}{2}}\phi (x)=\int_{-\infty }^{+\infty} K_k(x,y)\phi
(y)dy,\quad x\in \RR .
$$
So, in order to establish our result we must analyze the integral
$$
\int_0^\infty
t^{\frac{k}{2}-1}\left|\Big(\frac{d}{dx}+x\Big)^l(W_t(x,y))\right|dt,\quad
x,y \in \RR .
$$

A calculation shows that we can write, for every $l\in \NN$,
\begin{equation}\label{Dl}
\Big(\frac{d}{dx}+x\Big)^l=\sum_{0\leq \rho +\sigma \leq l}c_{\rho
,\sigma }^lx^\rho \frac{d^\sigma}{dx^\sigma}\;,
\end{equation}
where $c_{\rho ,\sigma }^l\in \RR$, $\rho ,\sigma \in \NN$, $0\leq
\rho +\sigma \leq l$ and $c_{0,l}^l=1$.

By making the change of variable $t=\log\frac{1+s}{1-s}$, we
obtain, for every $s\in (0,1) $ and $x,y \in \RR$,
$$
\Big(\frac{d}{dx}+x\Big)^lW_t(x,y)=\frac{1}{\sqrt{\pi
}}\left(\frac{1-s^2}{4s}\right)^{\frac{1}{2}}\sum_{0\leq \rho
+\sigma \leq l}c_{\rho ,\sigma }^lx^\rho
\frac{d^\sigma}{dx^\sigma}
\Big(e^{-\frac{1}{4}(s(x+y)^2+\frac{1}{s}(x-y)^2)}\Big),\quad l\in
\NN.
$$
Moreover, according to \cite[p. 50]{StTo2},
$\frac{d^\sigma}{dx^\sigma}
\Big(e^{-\frac{1}{4}(s(x+y)^2+\frac{1}{s}(x-y)^2)}\Big)$ is a
linear combination of terms of the form
$$
\left(s+\frac{1}{s}\right)^{b_1}e^{-\frac{1}{4}(s(x+y)^2+\frac{1}{s}(x-y)^2)}
\Big(s(x+y)+\frac{1}{s}(x-y)\Big)^{b_2},
$$
where $b_1,b_2\in \NN$ and $2b_1+b_2\leq \sigma$.

Assume that $l,\rho, \sigma ,b_1, b_2 \in \NN$, $0\leq l\leq k$,
$0\leq \rho +\sigma \leq l$ and $2b_1+b_2\leq \sigma$. Let us
consider
$$
\begin{array}{ll}
I_{l,\rho, \sigma}^{b_1,b_2}(x,y)&=\displaystyle x^\rho
\int_0^1\left(\log
\frac{1+s}{1-s}\right)^{\frac{k}{2}-1}\left(\frac{1-s^2}{s}\right)^{\frac{1}{2}}
\left(s+\frac{1}{s}\right)^{b_1}e^{-\frac{1}{4}(s(x+y)^2+\frac{1}{s}(x-y)^2)}\\
\\
&\displaystyle \quad \quad \times
\Big(s(x+y)+\frac{1}{s}(x-y)\Big)^{b_2}\frac{1}{1-s^2}ds,\quad
x,y\in \RR .
\end{array}
$$

We have that
$$
\begin{array}{ll}
|I_{l,\rho, \sigma}^{b_1,b_2}(x,y)|&\displaystyle \leq C
\left(\int_0^{\frac{1}{2}}s^{\frac{k}{2}-\frac{3}{2}-b_1-\frac{b_2}{2}-\frac{\rho
}{2}}e^{-c\frac{(x-y)^2}{s}}ds+\int_{\frac{1}{2}}^1(-\log
(1-s))^{\frac{k}{2}-1}\frac{1}{\sqrt{1-s}}ds\right)\\
\\
&\displaystyle \leq
C\left(\int_0^{\frac{1}{2}}s^{\frac{1}{2}(k-l-3)}e^{-c\frac{(x-y)^2}{s}}ds+1\right).
\end{array}
$$

Hence, according to \cite[Lemma 1.1]{StTo2},
$$
|I_{l,\rho, \sigma}^{b_1,b_2}(x,y)|\leq C\left\{\begin{array}{ll}
                                            1,&0\leq l\leq k-2,\\
                                            \displaystyle
                                            \frac{1}{\sqrt{|x-y|}},
                                            &l=k-1,\\
                                            \displaystyle
                                            \frac{1}{|x-y|},&l=k,
                                            \end{array}
\right.,\quad x,y\in \RR ,
$$
which allows us to conclude the result. \fin

We now complete the above result proving that, for every $\phi \in
C_c^\infty (\RR)$ and $k\in \NN$, $H^{-\frac{k}{2}}\phi$ is
$k$-times differentiable on $\RR$ and that
$\Big(\frac{d}{dx}+x\Big)^kH^{-\frac{k}{2}}\phi (x)$ is given by a
principal value integral, for every $x\in \RR$.

\begin{propo}\label{DHermite}
Let $\phi \in C_c^{\infty}(\RR )$ and $k \in \Bbb N$. Then
\begin{equation}\label{PV1}
\Big(\frac{d}{dx}+x\Big)^kH^{-\frac{k}{2}}\phi(x) = w_k\phi
(x)+\lim_{\varepsilon \rightarrow 0^+} \int_{|x-y|>\varepsilon}
R^{(k)}(x,y)\phi(y)\,dy,\;\; x\in \RR ,
\end{equation}
where $w_k=0$, if $k$ is odd, and $w_k=-2^{\frac{k}{2}}$, when $k$
is even.
\end{propo}

\noindent \underline{Proof}. We first observe that by Proposition
\ref{derivH} and since $\Big(\frac{d}{dx}+x\Big)^k=\sum_{{\rho
,\sigma \in \NN}\choose {0\leq \rho +\sigma \leq k}}c_{\rho
,\sigma}x^\rho \frac{d^\sigma}{dx^\sigma}$, for certain $c_{\rho
,\sigma}\in \RR$, $\rho ,\sigma \in \NN$, $0\leq \rho +\sigma \leq
k$, it is sufficient to prove that
\begin{equation}\label{PV1Dk}
\frac{d^k}{dx^k}H^{-\frac{k}{2}}\phi (x)= w_k\phi
(x)+\lim_{\varepsilon \rightarrow 0^+} \int_{|x-y|>\varepsilon}
\phi(y)\frac{1}{\Gamma (\frac{k}{2})}\int_0^\infty
t^{\frac{k}{2}-1}\frac{d^k}{dx^k}W_t(x,y)dtdy,\;\; x\in \RR ,
\end{equation}
where $w_k=0$, if $k$ is odd, and $w_k=-2^{\frac{k}{2}}$, when $k$
is even.

By making the change of variable $\displaystyle t=\log
\frac{1+s}{1-s}$ and by using (\ref{2.3}), we can write, for every
$s\in (0,1)$ and $x,y\in \RR$,
$$
\begin{array}{l}
\displaystyle \frac{d^k}{dx^k}W_{\log \frac{1+s}{1-s}}(x,y)=
\left(\frac{1-s^2}{4\pi s}\right)^{\frac{1}{2}}
\frac{d^k}{dx^k}\left[e^{-\frac{1}{4}(s(x+y)^2+\frac{1}{s}(x-y)^2)}\right]\vspace{3mm}\\
\displaystyle = \left(\frac{1-s^2}{4\pi s}\right)^{\frac{1}{2}}
\sum^k_{j=0} \dbinom{k}{j}
\frac{d^{j}}{dx^{j}}\left(e^{-\frac{s}{4}(x+y)^2}\right)
\frac{d^{k-j}}{dx^{k-j}}\left(e^{-\frac{1}{4s}(x-y)^2}\right)\vspace{3mm}\\
\displaystyle  =\left(\frac{1-s^2}{4\pi
s}\right)^{\frac{1}{2}}e^{-\frac{s}{4}(x+y)^2}
\frac{d^k}{dx^k}\left(e^{-\frac{1}{4s}(x-y)^2}\right)\vspace{3mm}\\
 + \displaystyle
\left(\frac{1-s^2}{4\pi s}\right)^{\frac{1}{2}}
\sum_{j=1}^k\dbinom{k}{j}\left(\sum_{l=0}^{E[\frac{j}{2}]}E_{j
,l}(x+y)^{j-2l}\Big(-\frac{s}{4}\Big)^{j
-l}e^{-\frac{s}{4}(x+y)^2}\right)\vspace{3mm}\\
\displaystyle \times \left(\sum_{m=0}^{E[\frac{k-j}{2}]}E_{k-j ,m}
(x-y)^{k-j-2m}\Big(-\frac{1}{4s}\Big)^{k-j
-m}e^{-\frac{1}{4s}(x-y)^2}\right).
\end{array}
$$

On the other hand, by proceeding as in the proof of Proposition
\ref{derivH} it follows that
\begin{equation}\label{PV1a}
\begin{array}{l}
\displaystyle \left|\int_{\frac{1}{2}}^1\left(\log
\frac{1+s}{1-s}\right)^{\frac{k}{2}-1}\frac{d^k}{dx^k}W_{\log
\frac{1+s}{1-s}}(x,y)\frac{2}{1-s^2}ds\right|\vspace{3mm}\\
\displaystyle \leq C\int_{\frac{1}{2}}^1(-\log
(1-s))^{\frac{k}{2}-1}\frac{1}{\sqrt{1-s}}ds\leq C,\quad x,y\in \RR
.
\end{array}
\end{equation}

We define, for every $j=1,...,k$, $0 \leq l \leq
E\Big[\frac{j}{2}\Big]$ and $0 \leq m \leq E\Big[\frac{k-j
}{2}\Big]$,
\begin{multline*}
F_{l,m}^j(x,y) = (x+y)^{j-2l}(x-y)^{k-j-2m} \\
\times \int_0^{\frac{1}{2}}
\left(\log\frac{1+s}{1-s}\right)^{\frac{k}{2}-1}\frac{e^{-\frac{1}{4}
(s(x+y)^2+\frac{1}{s}(x-y)^2)}}{s^{k-2j+l-m+\frac{1}{2}}}
\frac{1}{\sqrt{1-s^2}}ds,\;\;x,y\in \RR .
\end{multline*}

Since $\log\frac{1+s}{1-s}\sim s$, as $s\rightarrow 0^+$, and
$j\geq 1$, it follows that

$\begin{array}{ll} |F_{l,m}^j(x,y)| &\displaystyle \leq  C
|x+y|^{j-2l}|x-y|^{k-j-2m}\int_0^{\frac{1}{2}}
s^{-\frac{k}{2}-\frac{3}{2}+2j-l+m} e^{-\frac{1}{4}(s(x+y)^2
+\frac{1}{s}(x-y)^2)}\,ds \\
& \displaystyle \leq C \int_0^{\frac{1}{2}}
s^{j-\frac{3}{2}}e^{-c\frac{(x-y)^2}{s}}\,ds \leq C
\int_0^{\frac{1}{2}}s^{-\frac{1}{2}}\,ds \leq C,\quad x,y\in \RR .
\end{array}
$

Hence,
\begin{equation}\label{PV1b}
\begin{array}{l}
\displaystyle
\left|\int_0^{\frac{1}{2}}\left[\frac{d^k}{dx^k}W_{\log
\frac{1+s}{1-s}}(x,y)- \left(\frac{1-s^2}{4\pi
s}\right)^{\frac{1}{2}}e^{-\frac{s}{4}(x+y)^2}
\frac{d^k}{dx^k}\left(e^{-\frac{1}{4s}(x-y)^2}\right)\right]\right.\\
\\
\displaystyle \left.\quad \quad \times \left(\log
\frac{1+s}{1-s}\right)^{\frac{k}{2}-1}\frac{2}{1-s^2}ds\right|\leq
C, \quad x,y\in \RR .
\end{array}
\end{equation}

We now observe that mean value theorem leads to \vspace{3mm}

$\begin{array}{l} \displaystyle \left| e^{-\frac{s}{4}(x+y)^2}
\left(\log
\frac{1+s}{1-s}\right)^{\frac{k}{2}-1}\frac{1}{\sqrt{1-s^2}}-(2s)^{\frac{k}{2}-1}
\right|\leq \left|\left(\log
\frac{1+s}{1-s}\right)^{\frac{k}{2}-1}-(2s)^{\frac{k}{2}{-1}}\right|\vspace{3mm}\\
\hspace{2mm} +\displaystyle \left[\left|
\frac{1}{\sqrt{1-s^2}}-1\right|e^{-\frac{s}{4}(x+y)^2}+\left|
e^{-\frac{s}{4}(x+y)^2}-1\right|\right]\left|\left(\log
\frac{1+s}{1-s}\right)^{\frac{k}{2}-1}\right|\vspace{3mm}\\
\hspace{2mm} \displaystyle \leq
C\left(s^{\frac{k}{2}+1}+\Big(s^2e^{-\frac{s}{4}(x+y)^2}+(x+y)^2s\Big)\left|\left(\log
\frac{1+s}{1-s}\right)^{\frac{k}{2}-1}\right|\right),
\end{array}
$ \vspace{2mm}

\noindent for every $s\in (0,\frac{1}{2})$ and $x,y\in \RR$.

Then, since $\log \frac{1+s}{1-s}\sim s$, as $s\rightarrow 0^+$,
by using again (\ref{2.3}) we get
\begin{equation}\label{PV1c}
\begin{array}{l}
\left|\displaystyle
\int_0^{\frac{1}{2}}\left(\left(\frac{1-s^2}{4\pi
s}\right)^{\frac{1}{2}}e^{-\frac{s}{4}(x+y)^2}\left(\log\frac{1+s}{1-s}\right)^{\frac{k}{2}-1}
\frac{2}{1-s^2}-\frac{(2s)^{\frac{k}{2}-1}}{\sqrt{\pi
s}}\right)\frac{d^k}{dx^k}\left(e^{-\frac{(x-y)^2}{4s}}\right)
\,ds \right| \\\vspace{3mm} \quad \displaystyle \leq
C\int_0^{\frac{1}{2}} s^{\frac{k}{2}-\frac{1}{2}}(s+(x+y)^2)
\left|\frac{d^k}{dx^k}\left(e^{-\frac{(x-y)^2}{4s}}\right)\right|
\,ds \\\vspace{3mm} \displaystyle \quad \leq C
\sum_{n=0}^{E[\frac{k}{2}]} |x-y|^{k-2n} \int_0^{\frac{1}{2}}
s^{-\frac{k}{2}-\frac{1}{2}+n}(s+(x+y)^2)
e^{-\frac{(x-y)^2}{4s}}ds
\\\vspace{3mm} \displaystyle \quad  \leq C (1+(x+y)^2), \quad x,y
\in \RR.
\end{array}
\end{equation}

In view of properties (\ref{PV1a}), (\ref{PV1b}) and (\ref{PV1c}),
to establish the desired property (\ref{PV1Dk}) it is sufficient
to prove that
\begin{multline*}
\frac{d^k}{dx^k}\int_{-\infty}^{+\infty}
\phi(y)\frac{1}{\Gamma(\frac{k}{2})}\int_0^{\frac{1}{2}}
\frac{(2s)^{\frac{k}{2}-1}}{\sqrt{\pi s}}e^{-\frac{1}{4s}(x-y)^2}
\,dsdy \\
= w_k\phi (x)+\lim_{\varepsilon \rightarrow 0^+}
\int_{|x-y|>\varepsilon} \phi(y)
\frac{1}{\Gamma(\frac{k}{2})}\int_0^{\frac{1}{2}}\frac{(2s)^{\frac{k}{2}-1}}{\sqrt{\pi
s}}\frac{d^k}{dx^k}\left(e^{-\frac{1}{4s}(x-y)^2}\right) \,ds
dy,\, x \in \RR ,
\end{multline*}
where $w_k=0$, if $k$ is odd, and $w_k=-2^{\frac{k}{2}}$, when $k$
is even.

As earlier we can see that, for each $\phi \in C^\infty _c(\RR)$,
\begin{multline*}
\frac{d^{k-1}}{dx^{k-1}}\int_{-\infty }^{+\infty}
\phi(y)\int_0^{\frac{1}{2}}\frac{(2s)^{\frac{k}{2}-1}}{\sqrt{\pi
s}}e^{-\frac{1}{4s}(x-y)^2} \,dsdy \\
\quad = \int_{-\infty }^{+\infty} \phi(y)\int_0^{\frac{1}{2}}
\frac{(2s)^{\frac{k}{2}-1}}{\sqrt{\pi
s}}\frac{d^{k-1}}{dx^{k-1}}\left(e^{-\frac{1}{4s}(x-y)^2}\right)\,dsdy,\;\;
x\in \RR .
\end{multline*}
Let us represent by $\Phi$ the following function
$$
\Phi(x)=\frac{1}{\Gamma (\frac{k}{2})}\int_0^{\frac{1}{2}}
\frac{(2s)^{\frac{k}{2}-1}}{\sqrt{\pi
s}}\frac{d^{k-1}}{dx^{k-1}}\left(e^{-\frac{x^2}{4s}}\right)\,ds,\;\;
x\in \RR .
$$
By proceeding as above we can see that $\Phi \in L^1(\RR )$.
Indeed, (\ref{2.3}) leads to
$$
\Phi(x)=\frac{1}{\Gamma
(\frac{k}{2})}\sum_{l=0}^{E[\frac{k-1}{2}]}(-1)^
{k-1-l}E_{k-1,l}x^{k-1-2l}\int_0^{\frac{1}{2}}\frac{(2s)^{\frac{k}{2}-1}}{\sqrt{\pi
s}}\Big(\frac{1}{4s}\Big)^{k-1-l}e^{-\frac{x^2}{4s}}ds,\quad x\in
\RR .
$$
and then, according to \cite[Lemma 1.1]{StTo2},
\begin{equation}\label{Phi}
|\Phi(x)|\leq C\sum_{l=0}^{E[\frac{k-1}{2}]}|x|^{k-1-2l}
\int_0^{\frac{1}{2}}s^{-\frac{k}{2}+l-\frac{1}{2}}e^{-\frac{x^2}{4s}}ds\leq
Ce^{-\frac{x^2}{4}}\int_0^{\frac{1}{2}}\frac{e^{-c\frac{x^2}{s}}}{s}ds\leq
C\frac{e^{-\frac{x^2}{4}}}{\sqrt{|x|}},\quad x\in \RR .
\end{equation}
Hence $\Phi \in L^1(\RR )$. Moreover $\Phi \in C^\infty (\RR
\setminus \{0\})$.

Also, when $k$ is even, we can see that
\begin{equation}\label{limphi}
\lim_{\varepsilon \rightarrow 0}\Phi (\varepsilon
)=-2^{\frac{k}{2}-1}.
\end{equation}
In effect, if $k$ is even we can write
$$
\Phi(\varepsilon )=-\frac{1}{\Gamma (\frac{k}{2})\sqrt{\pi
}}\sum_{l=0}^{\frac{k}{2}-1}(-1)^lE_{k-1,l}\frac{\varepsilon ^
{k-1-2l}}{2^{\frac{3k}{2}-2l-1}}\int_0^{\frac{1}{2}}
\frac{e^{-\frac{\varepsilon ^2}{4s}}}{s^{\frac{k}{2}+
\frac{1}{2}-l}}ds,\quad \varepsilon \in \RR .
$$

Hence, the duplication formula (\cite[(1.2.3)]{Lebe}) allows us to
write
$$
\begin{array}{l}
\displaystyle  \lim_{\varepsilon \rightarrow 0}\Phi
(\varepsilon)=-\frac{1}{2^{\frac{k}{2}}\Gamma (\frac{k}{2})\sqrt{\pi
}}\lim_{\varepsilon \rightarrow
0}\sum_{l=0}^{\frac{k}{2}-1}(-1)^lE_{k-1,l}\int_{\frac{\varepsilon
^2}{2}}^\infty e^{-u}u^{\frac{k}{2}-
\frac{3}{2}-l}du\vspace{3mm}\\
\quad \displaystyle =\frac{-1}{2^{\frac{k}{2}}\Gamma
(\frac{k}{2})\sqrt{\pi
}}\sum_{l=0}^{\frac{k}{2}-1}(-1)^lE_{k-1,l}\Gamma
\Big(\frac{k-1}{2}-l\Big)=\frac{-(k-1)!}{2^{\frac{k}{2}-1}(\Gamma
(\frac{k}{2}))^2}\sum_{l=0}^{\frac{k}{2}-1}(-1)^l
{{\frac{k}{2}-1}\choose l}\frac{1}{k-1-2l}\vspace{3mm}\\
\quad \displaystyle =\frac{-(k-1)!}{2^{\frac{k}{2}-1}(\Gamma
(\frac{k}{2}))^2}\int_0^1(1-t^2)^{\frac{k}{2}-1}dt=
\frac{-(k-1)!}{2^{\frac{k}{2}}(\Gamma (\frac{k}{2}))^2}\frac{\Gamma
(\frac{k}{2})\Gamma (\frac{1}{2})}{\Gamma
(\frac{k+1}{2})}=-2^{\frac{k}{2}-1},
\end{array}
$$
and (\ref{limphi}) is thus established.

For every $x\in \RR$, we can write
$$
\begin{array}{l}
\displaystyle \frac{d}{dx}\int_{-\infty }^{+\infty } \phi (y)\Phi
(x-y)dy=\frac{d}{dx}\int_{-\infty }^{+\infty }\phi (x-y)\Phi
(y)dy=-\int_{-\infty }^{+\infty
}\frac{d}{dy}(\phi (x-y))\Phi (y)dy\vspace{3mm}\\
\hspace{2em}=\displaystyle -\lim _{\varepsilon \rightarrow
0^+}\int_{|y|>\varepsilon}\frac{d}{dy}(\phi (x-y))\Phi
(y)dy \vspace{3mm}\\
\hspace{2em}=\displaystyle \lim_{\varepsilon \rightarrow
0^+}\left[\int_{|y|>\varepsilon }\phi (x-y)\frac{d}{dy}\Phi
(y)dy-\phi (x+\varepsilon )\Phi (-\varepsilon )+\phi
(x-\varepsilon )\Phi
(\varepsilon )\right]\vspace{3mm}\\
\hspace{2em}=\displaystyle \lim_{\varepsilon \rightarrow
0^+}\left[\int_{|x-y|>\varepsilon }\phi (y)\Big(\frac{d}{dy}\Phi
\Big)(x-y)dy+\phi (x-\varepsilon )\Phi (\varepsilon )-\phi
(x+\varepsilon )\Phi (-\varepsilon )\right].
\end{array}
$$

Suppose now that $k$ is odd. Then $\Phi$ is an even function and
from (\ref{Phi}) we obtain, for every $x\in \RR$,
$$
\left|\phi (x-\varepsilon )\Phi (\varepsilon )-\phi (x+\varepsilon
)\Phi (-\varepsilon )\right|\leq C\varepsilon |\Phi (\varepsilon
)|\longrightarrow 0,\quad \mbox{ as }\varepsilon \rightarrow
0^+\;.
$$

On the other hand, assuming that $k$ is even, (\ref{limphi}) leads
to
$$
\lim_{\varepsilon \rightarrow 0^+}\phi (x-\varepsilon )\Phi
(\varepsilon )-\phi (x+\varepsilon )\Phi (-\varepsilon
)=\lim_{\varepsilon \rightarrow 0^+}(\phi (x+\varepsilon )+\phi
(x-\varepsilon ))\Phi (\varepsilon )=-2^{\frac{k}{2}}\phi (x),
$$
for every $x\in \RR$.

Hence,
$$
\displaystyle \frac{d}{dx}\int_{-\infty }^{+\infty}\phi (y)\Phi
(x-y)dy=w_k\phi (x)+\lim_{\varepsilon \rightarrow
0^+}\int_{|x-y|>\varepsilon } \phi (y)\frac{d}{dx}(\Phi
(x-y))dy,\quad x\in \RR ,
$$
where $w_k=0$, if $k$ is odd, and $w_k=-2^{\frac{k}{2}}$, when $k$
is even. Thus the proof is finished.

\fin

The following relation between the kernels $R_\alpha^{(k)}(x,y)$ and
$R^{(k)}(x,y)$, $x,y\in (0,\infty)$, is the key of our procedure in
order to establish that the $k$-order Riesz transform associated
with the Laguerre operator is a principal value integral operator.

\begin{propo}\label{kernels}
Let $\alpha>-1$ and $k\in \mathbb{N}$. We have that

(i) $\displaystyle |R_\alpha^{(k)}(x,y)|\leq
C\frac{y^{\alpha+\frac{1}{2}}}{x^{\alpha+\frac{3}{2}}}$,
$\displaystyle 0<y<\frac{x}{2}$.

(ii) $\displaystyle |R_\alpha^{(k)}(x,y)|\leq
C\frac{x^{\alpha+\frac{1}{2}}}{y^{\alpha+\frac{3}{2}}}$, $y>2x
\textup{ and } k \textup{ even, and}$  $\displaystyle
|R_\alpha^{(k)}(x,y)|\leq
C\frac{x^{\alpha+\frac{3}{2}}}{y^{\alpha+\frac{5}{2}}}$, $y>2x
\textup{ and } k \textup{ odd}$.

(iii) $\displaystyle \left|R_\alpha^{(k)}(x,y)- R^{(k)}(x,y)
\right|\leq C\frac {1}{x} \left(
1+\left(\frac{x}{|x-y|}\right)^{\frac{1}{2}}\right)$,
$\displaystyle \frac{x}{2}<y<2x$.
\end{propo}

\noindent \underline{Proof}. We first establish the following
formula that will be used later. For every $j\in \NN$, and
$t,x,y\in (0,\infty)$ we have that
\begin{equation}\label{DjI}
\begin{array}{l}
\displaystyle \frac{d^j}{dx^j}\left[ \left( \frac{2xy
e^{-t}}{1-e^{-2t}} \right)^{-\alpha} I_\alpha \left( \frac{2xy
e^{-t}}{1-e^{-2t}}\right) \right] \\
\\
\quad \displaystyle = \sum_{n=0}^{E[\frac{j}{2}]} E_{j,n}
\frac{x^{j-2n}}{2^{j-n}} \left( \frac{2y e^{-t}}{1-e^{-2t}}
\right)^{2(j-n)} \left( \frac{2xy e^{-t}}{1-e^{-2t}}
\right)^{-\alpha+n-j}I_{\alpha-n+j} \left( \frac{2xy
e^{-t}}{1-e^{-2t}}\right).
\end{array}
\end{equation}
Indeed, let $j\in \NN$. By using (\ref{2.3}) and since
$\frac{d}{dx} g(x) = \frac12 \left( \frac1u \frac{d}{du} \right)
[g(u^2)]_{|u=\sqrt{x}}$ we can write
$$
\begin{array}{l}
\displaystyle \frac{d^j}{dx^j}\left[ \left( \frac{2xy
e^{-t}}{1-e^{-2t}} \right)^{-\alpha} I_\alpha \left( \frac{2xy
e^{-t}}{1-e^{-2t}}\right) \right] \vspace{3mm}\\
\hspace{1cm}\displaystyle = \sum_{n=0}^{E[\frac{j}{2}]} E_{j,n}
x^{j-2n} \frac{d^{j-n}}{dz^{j-n}}\left( \left( \frac{2\sqrt{z}y
e^{-t}}{1-e^{-2t}} \right)^{-\alpha} I_\alpha \left(
\frac{2\sqrt{z}y e^{-t}}{1-e^{-2t}}\right)\right)_{\Big|z=x^2}\vspace{3mm}\\
\hspace{1cm}\displaystyle =\sum_{n=0}^{E[\frac{j}{2}]} E_{j,n}
\frac{x^{j-2n}}{2^{j-n}}
\Big(\frac{1}{x}\frac{d}{dx}\Big)^{j-n}\left( \left( \frac{2xy
e^{-t}}{1-e^{-2t}} \right)^{-\alpha} I_\alpha \left( \frac{2xy
e^{-t}}{1-e^{-2t}}\right)\right),\quad t,x,y\in (0,\infty ).
\end{array}
$$
Thus, (\ref{DjI}) can be easily deduced from property (P3).

Let us now prove $(i)$ and $(ii)$. Since
$\frac{d}{dx}+x=e^{-\frac{x^2}{2}}\frac{d}{dx}e^{\frac{x^2}{2}}$
we can write
$$
\begin{array}{ll}
\displaystyle \mathfrak{D}_\alpha^k W_t^\alpha (x,y)
&\displaystyle = x^{\alpha+\frac{1}{2}} e^{-\frac{x^2}{2}}
\frac{d^k}{dx^k}\Big(e^{\frac{x^2}{2}} x^{-\alpha-\frac{1}{2}}
W_t^\alpha (x,y)\Big)\vspace{3mm}\\
&\displaystyle =\left( \frac{2e^{-t}}{1-e^{-2t}}
\right)^{\frac{1}{2}} \left( \frac{2xy e^{-t}}{1-e^{-2t}}
\right)^{\alpha+\frac{1}{2}}
e^{-\frac{x^2}{2}-\frac{y^2}{2}\frac{1+e^{-2t}}{1-e^{-2t}}}\vspace{3mm}\\
&\displaystyle \times \frac{d^k}{dx^k}\left [\left( \frac{2xy
e^{-t}}{1-e^{-2t}} \right)^{-\alpha} I_\alpha \left( \frac{2xy
e^{-t}}{1-e^{-2t}}\right) e^{-x^2
\frac{e^{-2t}}{1-e^{-2t}}}\right],\quad t,x,y \in (0,\infty).
\end{array}
$$

By taking into account formulas (\ref{2.3}) and (\ref{DjI}) we get
$$
\begin{array}{l}
\displaystyle \frac{d^k}{dx^k} \left[\left( \frac{2xy
e^{-t}}{1-e^{-2t}} \right)^{-\alpha}I_\alpha \left( \frac{2xy
e^{-t}}{1-e^{-2t}}\right)e^{-x^2 \frac{e^{-2t}}{1-e^{-2t}}}\right]\vspace{3mm}\\
\displaystyle \;\;=\sum_{j=0}^k{k\choose j}\frac{d^j}{dx^j}
\left[\Big(\frac{2xye^{-t}}{1-e^{-2t}}\Big)^{-\alpha }I_\alpha
\Big(\frac{2xye^{-t}}{1-e^{-2t}}\Big)\right]\frac{d^{k-j}}{dx^{k-j}}
\left[e^{-x^2\frac{e^{-2t}}{1-e^{-2t}}}\right]\vspace{3mm}\\
 \;\;=\displaystyle
e^{-x^2\frac{e^{-2t}}{1-e^{-2t}}}\sum_{j=0}^k\sum_{n=0}^{E[\frac{j}{2}]}\sum_{m=0}^{
E[\frac{k-j}{2}]}{k\choose j}\frac{E_{j,n}E_{k-j,m}}{2^{j-n}}
\Big(\frac{2ye^{-t}}{1-e^{-2t}}\Big)^{2(j-n)}
\Big(\frac{-e^{-2t}}{1-e^{-2t}}\Big)^{k-j-m}\vspace{3mm}\\
\;\;\displaystyle \times
x^{k-2m-2n}\Big(\frac{2xye^{-t}}{1-e^{-2t}}\Big)^{-\alpha
-j+n}I_{\alpha +j-n}\Big(\frac{2xye^{-t}}{1-e^{-2t}}\Big),\quad
t,x,y\in (0,\infty ).
\end{array}
$$
Hence we obtain that
\begin{equation}\label{DW1}
\begin{array}{ll}
{\mathfrak{D}}_\alpha ^kW_t^\alpha (x,y)&=\displaystyle
\Big(\frac{2e^{-t}}{1-e^{-2t}}\Big)^{\frac{1}{2}}\Big(\frac{2xye^{-t}}{1-e^{-2t}}\Big)^{\alpha
+\frac{1}{2}}e^{-\frac{1}{2}(x^2+y^2)\frac{1+e^{-2t}}{1-e^{-2t}}}\\
&\\
&\displaystyle \times
\sum_{j=0}^k\sum_{n=0}^{E[\frac{j}{2}]}\sum_{m=0}^{E[\frac{k-j}{2}]}{k\choose
j}\frac{E_{j,n}E_{k-j,m}}{2^{j-n}}
\Big(\frac{2ye^{-t}}{1-e^{-2t}}\Big)^{2(j-n)}
\Big(\frac{-e^{-2t}}{1-e^{-2t}}\Big)^{k-j-m}\\
&\\
&\displaystyle \times
x^{k-2m-2n}\Big(\frac{2xye^{-t}}{1-e^{-2t}}\Big)^{-\alpha
-j+n}I_{\alpha +j-n}\Big(\frac{2xye^{-t}}{1-e^{-2t}}\Big),\quad
t,x,y\in (0,\infty ).
\end{array}
\end{equation}
By using property (P1) it follows that
$$
\begin{array}{l}
\displaystyle \left| \int_{0,\frac{2xye^{-t}}{1-e^{-2t}}\leq
1}^\infty t^{\frac{k}{2}-1}{\mathfrak{D}}_\alpha ^kW_t^\alpha
(x,y)dt\right|\leq C(xy)^{\alpha
+\frac{1}{2}}\sum_{j=0}^k\sum_{n=0}^{E[\frac{j}{2}]}\sum_{m=0}^{E[\frac{k-j}{2}]}
x^{k-2m-2n}y^{2(j-n)}\vspace{3mm}\\
\displaystyle \times \int_{0,\frac{2xye^{-t}}{1-e^{-2t}}\leq
1}^\infty
t^{\frac{k}{2}-1}e^{-\frac{1}{2}(x^2+y^2)\frac{1+e^{-2t}}{1-e^{-2t}}}
\Big(\frac{e^{-t}}{1-e^{-2t}}\Big)^{\alpha +1+j+k-2n-m}dt\vspace{3mm}\\
\leq \displaystyle C(xy)^{\alpha
+\frac{1}{2}}\sum_{j=0}^k\sum_{n=0}^{E[\frac{j}{2}]}\sum_{m=0}^{E[\frac{k-j}{2}]}
x^{k-2m-2n}y^{2(j-n)}\vspace{3mm}\\
\displaystyle \times \left(\int_0^1t^{-\frac{k}{2}-\alpha -2
-j+2n+m}e^{-c\frac{x^2+y^2}{t}}dt+e^{-c(x^2+y^2)}\int_1^\infty
t^{\frac{k}{2}-1}e^{-(\alpha +1)t}dt\right)
\end{array}
$$
Hence by taking into account \cite[Lemma 1.1]{StTo2} we conclude
that
\begin{equation}\label{IDW1ae}
\begin{array}{ll}
\left|\displaystyle \int_{0,\frac{2xye^{-t}}{1-e^{-2t}}\leq
1}^\infty t^{\frac{k}{2}-1}{\mathfrak{D}}_\alpha ^kW_t^\alpha
(x,y)dt\right|&\leq \displaystyle
C\sum_{j=0}^k\sum_{n=0}^{E[\frac{j}{2}]}\sum_{m=0}^{E[\frac{k-j}{2}]}
\frac{(xy)^{\alpha
+\frac{1}{2}}x^{k-2m-2n}y^{2(j-n)}}{(x^2+y^2)^{\frac{k}{2}+\alpha
+1+j-2n-m}}\vspace{3mm}\\
&\displaystyle \leq C\frac{(xy)^{\alpha
+\frac{1}{2}}}{(x^2+y^2)^{\alpha +1}} \leq
C\left\{\begin{array}{ll}
                 \displaystyle \frac{y^{\alpha +\frac{1}{2}}}{x^{\alpha +\frac{3}{2}}},&\displaystyle
                0<y<x\;,\\
                &\\
               \displaystyle \frac{x^{\alpha +\frac{1}{2}}}{y^{\alpha +\frac{3}{2}}},&y>x>0.
                \end{array}
\right.
\end{array}
\end{equation}
Note that if $k$ is odd we can improve the estimate when $y>x>0$
as follows
\begin{equation}\label{IDW1ao}
\left|\int_{0,\frac{2xye^{-t}}{1-e^{-2t}}\leq 1}^\infty
t^{\frac{k}{2}-1}{\mathfrak{D} }_\alpha ^kW_t^\alpha
(x,y)dt\right|\leq C
                \displaystyle \frac{(xy)^{\alpha +\frac{1}{2}}x}{(x^2+y^2)^{\alpha +\frac{3}{2}}}\leq
                \displaystyle \frac{x^{\alpha +\frac{3}{2}}}{y^{\alpha
                +\frac{5}{2}}}.
\end{equation}

Assume now that $\frac{2xye^{-t}}{1-e^{-2t}}\geq 1$. From
(\ref{DW1}) and property (P2) we get
$$
\begin{array}{ll}
\left|{\mathfrak{D}}_\alpha ^kW_t^\alpha (x,y)\right|&\leq
\displaystyle
C\sum_{j=0}^k\sum_{n=0}^{E[\frac{j}{2}]}\sum_{m=0}^{E[\frac{k-j}{2}]}
e^{-\frac{1}{2}(x^2+y^2)\frac{1+e^{-2t}}{1-e^{-2t}}+\frac{2xye^{-t}}{1-e^{-2t}}}\vspace{3mm}\\
&\displaystyle \times x^{k-2m-2n}y^{2(j-n)}
\Big(\frac{e^{-t}}{1-e^{-2t}}\Big)^{k-m-2n+j+\frac{1}{2}},\quad
t,x,y\in (0,\infty ).
\end{array}
$$
We also observe that
$$
-\frac{1}{2}(x^2+y^2)\frac{1+e^{-2t}}{1-e^{-2t}}+\frac{2xye^{-t}}{1-e^{-2t}}=
-\frac{(x-ye^{-t})^2+(y-xe^{-t})^2}{2(1-e^{-2t})}.
$$
Thus, if $0<y<\frac{x}{2}$, we can write
$$
\begin{array}{ll}
\left|{\mathfrak{D}}_\alpha ^kW_t^\alpha (x,y)\right|&\leq
\displaystyle
C\sum_{j=0}^k\sum_{n=0}^{E[\frac{j}{2}]}\sum_{m=0}^{E[\frac{k-j}{2}]}
e^{-\frac{x^2}{8(1-e^{-2t})}}x^{k-2m-2n+2(j-n)}\left(\frac{e^{-t}}{1-e^{-2t}}\right)^{k-m-2n+j+\frac{1}{2}}\vspace{3mm}\\
&\displaystyle \leq
Ce^{-c\frac{x^2}{1-e^{-2t}}}\left(\frac{e^{-t}}{1-e^{-2t}}\right)^{\frac{k}{2}+\frac{1}{2}},\quad
t\in (0,\infty ).
\end{array}
$$
Hence, if $-1<\alpha <-\frac{1}{2}$, \cite[Lemma 1.1]{StTo2} leads
to
$$
\begin{array}{l}
\displaystyle \left|\int_{0,\frac{2xye^{-t}}{1-e^{-2t}}\geq
1}^\infty t^{\frac{k}{2}-1}{\mathfrak{D} }_\alpha ^kW_t^\alpha
(x,y)dt\right|\leq
C\left(\int_0^1\frac{e^{-c\frac{x^2}{t}}}{t^{\frac{3}{2}}}dt+e^{-cx^2}\right)\leq
C\frac{1}{x}\leq C\frac{y^{\alpha +\frac{1}{2}}}{x^{\alpha
+\frac{3}{2}}},\;\;0<y<\frac{x}{2}.
\end{array}
$$
For $\alpha >-\frac{1}{2}$ we can proceed as follows.
$$
\begin{array}{l}
\displaystyle \left|\int_{0,\frac{2xye^{-t}}{1-e^{-2t}}\geq
1}^\infty t^{\frac{k}{2}-1}{\mathfrak{D} }_\alpha ^kW_t^\alpha
(x,y)dt\right|\leq C(xy)^{\alpha +\frac{1}{2}}\int_0^\infty
t^{\frac{k}{2}-1}e^{-c\frac{x^2}{1-e^{-2t}}}\left(\frac{e^{-t}}{1-e^{-2t}}\right)^{\frac{k}{2}+\alpha
+1}dt\vspace{3mm}\\
\displaystyle \leq C(xy)^{\alpha
+\frac{1}{2}}\left(\int_0^1\frac{e^{-c\frac{x^2}{t}}}{t^{\alpha
+2}}dt+e^{-cx^2}\right)\leq C\frac{(xy)^{\alpha
+\frac{1}{2}}}{x^{2\alpha +2}}\leq C\frac{y^{\alpha
+\frac{1}{2}}}{x^{\alpha +\frac{3}{2}}},\quad 0<y<\frac{x}{2}.
\end{array}
$$
In a similar way, if $0<2x<y$, we can write
$$
\begin{array}{l}
\displaystyle \left|\int_{0,\frac{2xye^{-t}}{1-e^{-2t}}\geq
1}^\infty t^{\frac{k}{2}-1}{\mathfrak{D} }_\alpha ^kW_t^\alpha
(x,y)dt\right|\leq C(xy)^{\alpha +\frac{3}{2}}\int_0^\infty
t^{\frac{k}{2}-1}e^{-c\frac{y^2}{1-e^{-2t}}}\left(\frac{e^{-t}}{1-e^{-2t}}\right)^{\frac{k}{2}+\alpha
+2}dt\vspace{3mm}\\
\displaystyle \leq C(xy)^{\alpha
+\frac{3}{2}}\left(\int_0^1\frac{e^{-c\frac{y^2}{t}}}{t^{\alpha
+3}}dt+e^{-cy^2}\right)\leq C\frac{(xy)^{\alpha
+\frac{3}{2}}}{y^{2\alpha +4}}\leq C\frac{x^{\alpha
+\frac{3}{2}}}{y^{\alpha +\frac{5}{2}}}.
\end{array}
$$
These estimations allow us to get
\begin{equation}\label{IDW1b}
\begin{array}{ll}
\left|\displaystyle \int_{0,\frac{2xye^{-t}}{1-e^{-2t}}\geq
1}^\infty t^{\frac{k}{2}-1}{\mathfrak{D}}_\alpha ^kW_t^\alpha
(x,y)dt\right|&\leq \displaystyle C\left\{\begin{array}{ll}
                 \displaystyle \frac{y^{\alpha +\frac{1}{2}}}{x^{\alpha +\frac{3}{2}}},&\displaystyle
                0<y<\frac{x}{2}\;,\\
                &\\
               \displaystyle \frac{x^{\alpha +\frac{3}{2}}}{y^{\alpha +\frac{5}{2}}},&y>2x>0.
                \end{array}
\right.
\end{array}
\end{equation}
Hence, by (\ref{IDW1ae}), (\ref{IDW1ao}) and (\ref{IDW1b}), $(i)$
and $(ii)$ are proved.

Next we establish statement $(iii)$. Observe first that, since
$\frac{d}{dx}+x=e^{-\frac{x^2}{2}}\frac{d}{dx}e^{\frac{x^2}{2}}$,
$$
\begin{array}{l}
\displaystyle \mathfrak{D}_\alpha^k W_t^\alpha (x,y)= \displaystyle
x^{\alpha +\frac{1}{2}}\Big(\frac{d}{dx}+x\Big)^k[x^{-\alpha
-\frac{1}{2}}W_t^\alpha (x,y)]\vspace{3mm}\\
\displaystyle =\sqrt{2\pi }x^{\alpha
+\frac{1}{2}}\Big(\frac{d}{dx}+x\Big)^k\left[x^{-\alpha
-\frac{1}{2}}e^{\frac{- 2xye^{-t}}{1-e^{-2t}}}\left( \frac{2xy
e^{-t}}{1-e^{-2t}}\right)^{\frac{1}{2}} I_\alpha \left( \frac{2xy
e^{-t}}{1-e^{-2t}} \right) W_t(x,y)\right]\vspace{3mm}\\
\displaystyle =\sqrt{2\pi }\left( \frac{2xy
e^{-t}}{1-e^{-2t}}\right)^{\alpha +\frac{1}{2}}\sum_{j=0}^k{k\choose
j}\frac{d^j}{dx^j}\left[e^{\frac{- 2xye^{-t}}{1-e^{-2t}}}\left(
\frac{2xy e^{-t}}{1-e^{-2t}}\right)^{-\alpha } I_\alpha \left(
\frac{2xy
e^{-t}}{1-e^{-2t}} \right)\right]\Big(\frac{d}{dx}+x\Big)^{k-j}W_t(x,y)\vspace{3mm}\\
\displaystyle = \sqrt{2\pi }\left( \frac{2xy
e^{-t}}{1-e^{-2t}}\right)^{\alpha
+\frac{1}{2}}\sum_{j=0}^k{k\choose
j}\Big(\frac{d}{dx}+x\Big)^{k-j}W_t(x,y)\vspace{3mm}\\
\displaystyle \times \sum_{l=0}^j(-1)^{j-l}{j\choose l}\left(\frac{
2ye^{-t}}{1-e^{-2t}}\right)^{j-l}e^{-\frac{
2xye^{-t}}{1-e^{-2t}}}\frac{d^l}{dx^l}\left(\left( \frac{2xy
e^{-t}}{1-e^{-2t}}\right)^{-\alpha } I_\alpha \left( \frac{2xy
e^{-t}}{1-e^{-2t}} \right)\right), \;\:t,x,y\in (0,\infty ).
\end{array}
$$
Hence, by using formula (\ref{DjI}) we obtain that, for every
$t,x,y\in (0,\infty )$,
\begin{equation}\label{DW2}
\begin{array}{l}
\displaystyle\displaystyle \mathfrak{D}_\alpha^k W_t^\alpha (x,y)
=\sqrt{2\pi} e^{-\frac{2xye^{-t}}{1-e^{-2t}}}\sum_{j=0}^k(-1)^j
{k\choose j} \Big(\frac{d}{dx}+x\Big)^{k-j} (W_t(x,y)) \left(
\frac{2y
e^{-t}}{1-e^{-2t}} \right)^j\\
\\
\displaystyle \times \sum_{n=0}^{E[\frac{j}{2}]}
\sum_{l=2n}^j(-1)^l{j \choose l}  \frac{E_{l,n}}{2^{l-n}} \left(
\frac{2xy e^{-t}}{1-e^{-2t}} \right)^{-n}\left( \frac{2xy
e^{-t}}{1-e^{-2t}} \right)^{\frac{1}{2}} I_{\alpha-n+l} \left(
\frac{2xy e^{-t}}{1-e^{-2t}} \right).
\end{array}
\end{equation}

Let us consider now $x,y,t\in (0,\infty )$ such that
$\frac{2xye^{-t}}{1-e^{-2t}}\geq 1$. By taking into account
property (P2) and (\ref{DW2}) we can write
$$
\begin{array}{l}
\displaystyle {\mathfrak{D}}_\alpha ^kW_t^\alpha (x,y)
=\Big(\frac{d}{dx}+x\Big)^k(W_t(x,y))
\left(1+O\left(\frac{1-e^{-2t}}{xye^{-t}}\right)\right)\vspace{3mm}\\
\displaystyle +\sum _{j=1}^k(-1)^j{k\choose
j}\Big(\frac{d}{dx}+x\Big)^{k-j}(W_t(x,y))\left(\frac{2ye^{-t}}{1-e^{-2t}}\right)^j
\sum_{n=0}^{E[\frac{j}{2}]} \sum_{l=2n}^j(-1)^l{j \choose
l}\frac{E_{l,n}}{2^{l-n}} \left( \frac{1-e^{-2t}}{2xy e^{-t}}
\right)^n\vspace{3mm}\\
\displaystyle \times
\left(\sum_{r=0}^{E[\frac{j}{2}]}\frac{(-1)^r[\alpha
+l-n,r]}{2^r}\left(\frac{1-e^{-2t}}{2xy e^{-t}} \right)^r +O\left(
\left(\frac{1-e^{-2t}}{xy e^{-t}}
\right)^{E[\frac{j}{2}]+1}\right)\right)\vspace{3mm}\\
=\displaystyle \Big(\frac{d}{dx}+x\Big)^k(W_t(x,y))+\sum
_{j=1}^k(-1)^j{k\choose
j}\Big(\frac{d}{dx}+x\Big)^{k-j}(W_t(x,y))\left(\frac{2ye^{-t}}{1-e^{-2t}}\right)^j\vspace{3mm}\\
\displaystyle \times \sum_{n=0}^{E[\frac{j}{2}]}
\sum_{l=2n}^j\sum_{r=0}^{E[\frac{j}{2}]}(-1)^{l+r}{j \choose
l}\frac{E_{l,n}}{2^{l-n}}\frac{[\alpha
+l-n,r]}{2^r}\left(\frac{1-e^{-2t}}{2xy e^{-t}}
\right)^{n+r}\vspace{3mm}\\
\displaystyle +\sum
_{j=0}^k(-1)^j{k\choose
j}\Big(\frac{d}{dx}+x\Big)^{k-j}(W_t(x,y))O\left( \left(\frac{y
e^{-t}}{1-e^{-2t}}\right)^{j-E[\frac{j}{2}]-1}\frac{1}{x^{E[\frac{j}{2}]+1}}\right).\vspace{3mm}\\
\end{array}
$$
Lemma \ref{technical} allows us to see that, for every $j\in \NN$,
$j=1,...,k$,
$$
\begin{array}{l}
\displaystyle \sum_{n=0}^{E[\frac{j}{2}]}
\sum_{l=2n}^j\sum_{r=0}^{E[\frac{j}{2}]}(-1)^{l+r}{j \choose
l}\frac{E_{l,n}}{2^{l-n}}\frac{[\alpha
+l-n,r]}{2^r}\left(\frac{1-e^{-2t}}{2xy e^{-t}} \right)^{n+r}\vspace{3mm}\\
\displaystyle = \sum_{n=0}^{E[\frac{j}{2}]}
\sum_{l=2n}^j\sum_{m=n}^{E[\frac{j}{2}]+n}(-1)^{l+m-n}{j \choose
l}\frac{E_{l,n}}{2^{l-n}}\frac{[\alpha
+l-n,m-n]}{2^{m-n}}\left(\frac{1-e^{-2t}}{2xy e^{-t}} \right)^m\vspace{3mm}\\
\displaystyle
=\sum_{m=0}^{E[\frac{j}{2}]}\left(\frac{1-e^{-2t}}{4xy e^{-t}}
\right)^m\sum_{n=0}^m\sum_{l=2n}^j(-1)^{l+m-n}{j \choose
l}\frac{E_{l,n}}{2^{l-2n}}[\alpha
+l-n,m-n]\vspace{3mm}\\
\displaystyle
+\sum_{m=E[\frac{j}{2}]+1}^{2E[\frac{j}{2}]}\left(\frac{1-e^{-2t}}{4xy
e^{-t}} \right)^m\sum_{n=m-E[\frac{j}{2}]}^{E[\frac{j}{2}]}
\sum_{l=2n}^j(-1)^{l+m-n}{j \choose
l}\frac{E_{l,n}}{2^{l-2n}}[\alpha
+l-n,m-n]\vspace{3mm}\\
\displaystyle =O\left(\left(\frac{1-e^{-2t}
}{xye^{-t}}\right)^{E[\frac{j}{2}]+1}\right).
\end{array}
$$

Hence, it follows that
$$
\begin{array}{ll}
\displaystyle {\mathfrak{D}}_\alpha ^kW_t^\alpha
(x,y)&\displaystyle =\Big(\frac{d}{dx}+x\Big)^kW_t(x,y)\\
&\\
&\displaystyle +\sum _{j=0}^k(-1)^j{k\choose
j}\Big(\frac{d}{dx}+x\Big)^{k-j}(W_t(x,y))O\left( \left(\frac{y
e^{-t}}{1-e^{-2t}}\right)^{j-E[\frac{j}{2}]-1}\frac{1}{x^{E[\frac{j}{2}]+1}}\right).
\end{array}
$$
Assume that $0<\frac{x}{2}<y<2x$. In order to establish $(iii)$ we
now proceed as in the proof of Proposition \ref{derivH}. First note
that by formula (\ref{Dl})
$$
\begin{array}{l}
\displaystyle \left|{\mathfrak{D}}_\alpha ^kW_t^\alpha
(x,y)-\Big(\frac{d}{dx}+x\Big)^kW_t(x,y)\right|\vspace{3mm}\\
\quad \displaystyle \leq C\sum_{j=0}^k\sum _{0\leq \rho +\sigma \leq
k-j}x^\rho \left|\frac{d^\sigma}{dx^\sigma}
W_t(x,y)\right|\left(\frac{ye^{-t}}{1-e^{-2t}}\right)^{j-E[\frac{j}{2}]-1}\frac{1}{x^{E[\frac{j}{2}]+1}}\;.
\end{array}
$$

Assume that $j,\rho, \sigma , b_1,b_2\in \NN$, $0\leq j\leq k$,
$0\leq \rho +\sigma \leq k-j$ and $2b_1+b_2\leq \sigma$. According
to \cite[p. 50]{StTo2} and by making the change of variable
$t=\log \frac{1+s}{1-s}$, we must analyze the following integral.
$$
\begin{array}{ll}
I_{\rho ,\sigma ,j}^{b_1,b_2}(x,y)&\displaystyle =\frac{x^\rho
y^j}{(xy)^{ 1+E[\frac{j}{2}]}}\int_{0,\frac{(1-s^2)xy}{2s}\geq 1}^1
\left(\log
\frac{1+s}{1-s}\right)^{\frac{k}{2}-1}\left(\frac{1-s^2}{s}\right)^{j-E[\frac{j}{2}]-\frac{1}{2}}
\vspace{3mm}\\
&\times \displaystyle
\left(s+\frac{1}{s}\right)^{b_1}e^{-\frac{1}{4}(s(x+y)^2+\frac{1}{s}(x-y)^2)}
\left(s(x+y)+\frac{1}{s}(x-y)\right)^{b_2}\frac{ds}{1-s^2}\vspace{3mm}\\
&=J_{\rho ,\sigma ,j}^{b_1,b_2}(x,y)+H_{\rho ,\sigma
,j}^{b_1,b_2}(x,y),\quad x,y\in (0,\infty ),
\end{array}
$$
where $J$ and $H$ are defined as $I$ but replacing the integral over
$(0,1)$ by the integral over $\Big(0,\frac{1}{2}\Big)$ and
$\Big(\frac{1}{2},1\Big)$, respectively.

Since $\log \frac{1+s}{1-s}\sim s$, as $s\rightarrow 0^+$, it
follows that
$$
\begin{array}{ll}
J_{\rho ,\sigma ,j}^{b_1,b_2}(x,y)&\displaystyle \leq C \frac{x^\rho
y^j}{(xy)^{ 1+E[\frac{j}{2}]}}\vspace{3mm}\\
&\times \displaystyle \int_{0,\frac{(1-s^2)xy}{2s}\geq
1}^{\frac{1}{2}}s^{\frac{k}{2}-\frac{1}{2}-j+E[\frac{j}{2}]-b_1}
e^{-\frac{1}{4}(s(x+y)^2+\frac{(x-y)^2}{s})}\left|s(x+y)+\frac{(x-y)}{s}\right|^{b_2}ds\vspace{3mm}\\
&\displaystyle \leq C\frac{y^j}{(xy)^{
1+E[\frac{j}{2}]}}\int_{0,\frac{(1-s^2)xy}{2s}\geq
1}^{\frac{1}{2}}s^{\frac{1}{2}(k-j-2b_1-b_2-\rho
)+E[\frac{j}{2}]-\frac{j}{2}-\frac{1}{2}}
e^{-c\frac{(x-y)^2}{s}}ds\vspace{3mm}\\
\end{array}
$$

$$
\begin{array}{ll}
&\displaystyle \leq C \frac{y^j}{(xy)^{
1+E[\frac{j}{2}]}}\int_{0,\frac{(1-s^2)xy}{2s}\geq
1}^{\frac{1}{2}}s^{E[\frac{j}{2}]-\frac{j}{2}-\frac{1}{2}}e^{-c\frac{(x-y)^2}{s}}ds.
\end{array}
$$

By taking into account that $0<\frac{x}{2}<y<2x$ and using
\cite[Lemma 1.1]{StTo2} we get
$$
\begin{array}{ll}
J_{\rho ,\sigma ,j}^{b_1,b_2}(x,y)&\displaystyle \leq C
\frac{x^{j-2E[\frac{j}{2}]-\frac{3}{2}}}{\sqrt{x}}\int_{0,\frac{(1-s^2)xy}{2s}\geq
1}^{\frac{1}{2}}s^{E[\frac{j}{2}]-\frac{j}{2}-\frac{1}{2}}e^{-c\frac{(x-y)^2}{s}}ds\vspace{3mm}\\
&\displaystyle \leq C
\frac{1}{\sqrt{x}}\int_0^{\frac{1}{2}}\frac{e^{-c\frac{(x-y)^2}{s}}}{s^{\frac{5}{4}}}ds
\leq C\frac{1}{x}\left(\frac{x}{|x-y|}\right)^{\frac{1}{2}}\;.
\end{array}
$$

On the other hand, since that $\log \frac{1+s}{1-s}\sim
-\log(1-s)$, as $s\rightarrow 1^-$, we have that
$$
\begin{array}{ll}
H_{\rho ,\sigma ,j}^{b_1,b_2}(x,y)&\displaystyle \leq C \frac{x^\rho
y^j}{(xy)^{
1+E[\frac{j}{2}]}}\int_{\frac{1}{2},\frac{(1-s^2)xy}{2s}\geq 1}^1
(-\log
(1-s))^{\frac{k}{2}-1}(1-s)^{j-E[\frac{j}{2}]-\frac{3}{2}}e^{-cs(x+y)^2}ds\vspace{3mm}\\
&\displaystyle \leq Ce^{-c(x+y)^2}\int_{\frac{1}{2}}^1(-\log
(1-s))^{\frac{k}{2}-1}(1-s)^{j-\frac{1}{2}}ds\leq
Ce^{-c(x+y)^2},\quad x,y\in (0,\infty ).
\end{array}
$$

Hence we conclude that, if $0<\frac{x}{2}<y<2x$,
\begin{equation}\label{IDW2a}
\left|\displaystyle \int_{0,\frac{2xye^{-t}}{1-e^{-2t}}\geq
1}^\infty t^{\frac{k}{2}-1}\Big({\mathfrak{D}}_\alpha ^kW_t^\alpha
(x,y)-\Big(\frac{d}{dx}+x\Big)^kW_t(x,y)\Big)dt\right|\leq C
\frac{1}{x}\left(\frac{x}{|x-y|}\right)^{\frac{1}{2}}\;.
\end{equation}

Also, by using again (\ref{2.3}) we obtain, for each $t,x,y\in
(0,\infty )$,
$$
\begin{array}{ll}
\displaystyle \Big(\frac{d}{dx}+x\Big)^kW_t(x,y) &\displaystyle
=e^{-\frac{x^2}{2}}\frac{d^k}{dx^k}\left[e^{\frac{x^2}{2}}W_t(x,y)\right]
\vspace{3mm}\\
&\displaystyle =W_t(x,y)
\sum_{j=0}^k\sum_{l=0}^{E[\frac{j}{2}]}{k\choose
j}E_{j,l}x^{j-2l}\Big(\frac{2ye^{-t}}{1-e^{-2t}}\Big)^{k-j}
\Big(\frac{-e^{-2t}}{1-e^{-2t}}\Big)^{j-l}.
\end{array}
$$

Hence it follows that, when $0<\frac{x}{2}<y<2x$,
\begin{equation}\label{IDW1H}
\begin{array}{l}
\displaystyle \left|\int_{0,\frac{2xye^{-t}}{1-e^{-2t}}\leq
1}^\infty t^{\frac{k}{2}-1}\Big(\frac{d}{dx}+x\Big)^kW_t(x,y)dt\right|\vspace{3mm}\\
\displaystyle \leq C \sum_{l=0}^{E[\frac{k}{2}]}x^{k-2l}
\left(\int_0^1t^{-\frac{k}{2}-\frac{3}{2}+l}e^{-c\frac{x^2}{t}}dt
+e^{-cx^2}\int_1^\infty
t^{\frac{k}{2}-1}e^{-\frac{t}{2}}dt\right)\leq C\frac{1}{x}.
\end{array}
\end{equation}

The estimations (\ref{IDW1ae}), (\ref{IDW2a}) and (\ref{IDW1H})
allow us to finish the proof of $(iii)$. \fin

By proceeding as above and having in mind Proposition \ref{derivH}
we can see that, for every $k \in \Bbb N$ and $\phi \in
C_c^{\infty}(0,\infty )$, the function $L^{-\frac{k}{2}}\phi$ is
$(k-1)$-times differentiable on $(0,\infty)$ and $k$-times
differentiable on $(0,\infty)\setminus \mbox{supp}\;\phi$. Moreover,
\begin{equation}\label{derivL}
\mathfrak{D}_{\alpha}^{\ell}L^{-\frac{k}{2}}\phi (x)=\frac{1}{\Gamma
(\frac{k}{2})}\int_0^{\infty} \phi(y)\int_0^\infty
t^{\frac{k}{2}-1}\mathfrak{D}_{\alpha}^{\ell}W_t^\alpha (x,y)dtdy,
\end{equation}
for every $x \in (0,\infty)$, when $\ell = 0,1,\ldots, k-1$, and for
every $x\in (0,\infty) \setminus \mbox{supp}\;\phi$, when $\ell =
k$.

We now prove that, for every $\phi \in C_c^{\infty}(0,\infty )$
and $k \in \Bbb N$, $L^{-\frac{k}{2}}_{\alpha} \phi$ is $k$-times
differentiable on $(0,\infty)$ and that $\mathfrak{D}_{\alpha}^k
L^{-\frac{k}{2}}\phi$ is a principal value integral operator
(modulus a constant times of the function when $k$ is even).

\begin{propo}\label{DL}
Let $\alpha > -1$, $k \in \Bbb N$ and $\phi \in
C_c^{\infty}(0,\infty)$. Then
$$
\mathfrak{D}_{\alpha}^k L^{-\frac{k}{2}}_{\alpha}\phi (x) =
w_k\phi (x)+\lim_{\varepsilon \rightarrow 0^+}
\int_{0,|x-y|>\varepsilon}^\infty
 R_\alpha ^{(k)}(x,y)\phi(y) \,dy, \;\;\, x
\in (0,\infty),
$$
where $w_k=0$, when $k$ is odd and $w_k=-2^{\frac{k}{2}}$, when
$k$ is even.
\end{propo}

Proof. For every $x\in (0,\infty )$, Proposition \ref{kernels}
implies that
$$
\begin{array}{l}
\displaystyle \lim_{\varepsilon \rightarrow
0^+}\int_{0,|x-y|>\varepsilon}^\infty R_\alpha ^{(k)}(x,y)\phi
(y)dy\vspace{3mm}\\
\displaystyle =\lim_{\varepsilon \rightarrow 0^+}
\int_{0,|x-y|>\varepsilon}^\infty (R_\alpha
^{(k)}(x,y)-R^{(k)}(x,y))\phi (y)dy+\lim_{\varepsilon \rightarrow
0^+}\int_{0,|x-y|>\varepsilon}^\infty R^{(k)}(x,y)\phi
(y)dy\vspace{3mm}\\
\displaystyle =\int_0^\infty \left(\mathfrak{D}_\alpha ^kK_{\alpha
,k}(x,y)-\Big(\frac{d}{dx}+x\Big)^kK_k(x,y)\right)\phi (y)dy
+\lim_{\varepsilon \rightarrow
0^+}\int_{0,|x-y|>\varepsilon}^\infty
R^{(k)}(x,y)\phi (y)dy\vspace{3mm}\\
\displaystyle =\frac{d}{dx}\left(\int_0^\infty
\left[\mathfrak{D}_\alpha ^{k-1}K_{\alpha
,k}(x,y)-\Big(\frac{d}{dx}+x\Big)^{k-1}K_k(x,y)\right]\phi (y)dy\right)\vspace{3mm}\\
\displaystyle +\left(x-\frac{\alpha
+\frac{1}{2}}{x}\right)\int_0^\infty \mathfrak{D}_\alpha
^{k-1}(K_{\alpha ,k}(x,y))\phi (y)dy-x\int_0^\infty
\Big(\frac{d}{dx}+x\Big)^{k-1}(K_k(x,y))\phi (y)dy\vspace{3mm}\\
\displaystyle +\lim_{\varepsilon \rightarrow
0^+}\int_{0,|x-y|>\varepsilon}^\infty R^{(k)}(x,y)\phi (y)dy.
\end{array}
$$

By taking into account Propositions \ref{derivH} and
\ref{DHermite} and (\ref{derivL}) we can conclude that
$$
\begin{array}{l}
\displaystyle \lim_{\varepsilon \rightarrow
0^+}\int_{0,|x-y|>\varepsilon}^\infty R_\alpha ^{(k)}(x,y)\phi
(y)dy\vspace{3mm}\\
\displaystyle \quad \quad =\frac{d}{dx}\left(\int_0^\infty
\left[\mathfrak{D}_\alpha ^{k-1}K_{\alpha
,k}(x,y)-\Big(\frac{d}{dx}+x\Big)^{k-1}K_k(x,y)\right]\phi (y)dy\right)\vspace{3mm}\\
\displaystyle  \quad \quad +\left(x-\frac{\alpha
+\frac{1}{2}}{x}\right)\int_0^\infty \mathfrak{D}_\alpha
^{k-1}K_{\alpha ,k}(x,y)\phi (y)dy-w_k\phi (x)
\end{array}
$$

$$
\begin{array}{l}
\displaystyle \quad \quad +\frac{d}{dx}\int_0^\infty
\Big(\frac{d}{dx}+x\Big)^{k-1}K_k(x,y)\phi (y)
dy=\mathfrak{D}_\alpha ^kL_\alpha ^{-\frac{k}{2}}\phi (x)-w_k\phi
(x),\quad x\in (0,\infty ),
\end{array}
$$
where $w_k=0$, for $k$ odd, and $w_k=-2^{\frac{k}{2}}$, when $k$
is even. Thus the proof is finished. \fin

We now prove the main result of the paper.

\noindent {\it Proof of Theorem \ref{main}}. We consider the maximal
operator associated with $R_\alpha ^{(k)}$ defined by
$$
R_{\alpha ,*}^{(k)}f(x)=\sup_{\varepsilon
>0}\left|\int_{0,|x-y|>\varepsilon }^\infty R_\alpha ^{(k)}(x,y)
f(y)dy\right|,\quad f\in C^\infty _c(0,\infty ), x\in (0,\infty ).
$$

According to Proposition \ref{kernels} we get
$$
R_{\alpha ,*}^{(k)}f(x)\leq C(H_0^{\alpha
+\frac{1}{2}}(|f|)(x)+H_\infty ^{\alpha +\frac{1}{2}+\delta
_k}(|f|)(x)+R_{{\rm loc},*}^{(k)}(f)(x)+N(f)(x)),
$$
where $\delta _k=1$, when $k$ is odd, $\delta _k=0$, when $k$ is
even,
$$
R_{{\rm loc},*}^{(k)}(f)(x)=\sup_{\varepsilon
>0}\left|\int_{\frac{x}{2},|x-y|>\varepsilon }^{2x}\Big(\frac{d}{dx}+x\Big)^k
K_k(x,y)f(y)dy\right|,
$$
and
$$
N(f)(x)=\int_{\frac{x}{2}}^{2x}f(y)\frac{1}{y}\left(1+\Big(\frac{x}{|x-y|}
\Big)^{\frac{1}{2}}\right)dy.
$$

By \cite[Lemma 3.1]{ChiPola} $H_0^{\alpha +\frac{1}{2}}$ is of
strong type $(p,p)$ with respect to $x^\delta dx$, when
$1<p<\infty$ and $\delta <\Big(\alpha +\frac{3}{2}\Big)p-1$, and
of weak type $(1,1)$ when $\delta \leq \alpha +\frac{1}{2}$. Also
from \cite[Lemma 3.2]{ChiPola} the operator $H_\infty ^{\alpha
+\frac{1}{2}+\delta _k}$ is of strong type $(p,p)$ for $x^\delta
dx$, when $1<p<\infty$ and $-\Big(\alpha
+\frac{1}{2}\Big)p-1<\delta$, and of weak type $(1,1)$ with
respect to $x^\delta dx$ when $-\alpha -\frac{5}{2}\leq \delta$,
if $k$ is odd; and, in the case that $k$ is even, when $\delta
\geq -\alpha -\frac{3}{2}$, and $\alpha \not=-\frac{1}{2}$ and
when $\delta >-1$ and $\alpha =-\frac{1}{2}$.

On the other hand, by using Jensen inequality we can see that the
operator $N$ is bounded from $L^p((0,\infty ),x^\delta dx)$ into
itself, for every $1\leq p<\infty$ and $\delta \in \RR$.

In \cite{StTo2} it was established that the kernel $R^{(k)}(x,y)$,
$x,y\in \RR$, is a Calder\'on-Zygmund kernel. Then, according to
\cite[Theorem 4.3]{NoSt1}, the operator $R_{{\rm loc},*}^{(k)}$ is
of strong type $(p,p)$, $1<p<\infty$, and of weak type $(1,1)$ with
respect to $x^\delta dx$, for every $\delta \in \RR$.

Then we conclude that $R_{\alpha ,*}^{(k)}$ defines an operator of
strong type $(p,p)$ for $x^\delta dx$ when $1<p<\infty$ and
$-\Big(\alpha +\frac{1}{2}+\delta _k\Big)p-1<\delta<\Big(\alpha
+\frac{3}{2}\Big)p-1$. We have also that $R_\alpha ^{(k)}$ is of
weak type $(1,1)$ for $x^\delta dx$ when $-\alpha -\frac{5}{2}\leq
\delta \leq \alpha +\frac{1}{2}$, if $k$ is odd. When $k$ is even
the maximal operator $R_{\alpha ,*}^{(k)}$ is of weak type (1,1)
with respect to $x^\delta dx$, for $-\alpha -\frac{3}{2}\leq
\delta \leq \alpha +\frac{1}{2}$ and $\alpha \not=-\frac{1}{2}$,
and for $-1<\delta \leq 0$, when $\alpha =-\frac{1}{2}$.

By using Proposition \ref{DL} and density arguments we can conclude
the proof of this theorem in a standard way. \fin


\begin{thebibliography}{99}

\bibitem{AMST} {I. Abu-Falahah, R. Mac\'{\i}as, C. Segovia and J.L.
Torrea}, Transferring strong boundedness among Laguerre orthogonal
systems, {\em preprint 2007}.

\bibitem{AnMu} { K.F. Andersen and B. Muckenhoupt}, Weighted weak type Hardy inequalities
with applications to Hilbert transforms and maximal functions, {\em
Studia Math.} 72 (1982), 9--26.

\bibitem{BFMR} {J.J. Betancor, J.C. Fari\~na, T. Mart\'{\i}nez and
L. Rodr\'{\i}guez-Mesa}, Higher order Riesz transforms associated with
Bessel operators, to appear in {\em Arkiv f\"or Matematik}.

\bibitem{BFRST} {J.J. Betancor, J.C. Fari\~na, L. Rodr\'{\i}guez-Mesa, A. Sanabria-Garc\'{\i}a and J.L.
Torrea}, Transference between Laguerre and Hermite settings, {\em
J. Funct. Anal.} 254 (2008), 826-850.

\bibitem{ChiPola} {A. Chicco Ruiz and E. Harboure}, Weighted norm
inequalities for the heat-diffusion Laguerre's semigroups, {\em
Math. Z.} 257 (2007), 329--354.

\bibitem{ChiPola2} {A. Chicco Ruiz and E. Harboure}, Boundedness
with weights for fractional integrals associated with Laguerre
expansions, preprint.

\bibitem{FaGuSc} {E. Fabes, C. Guti\'errez and R. Scotto},
Weak-type estimates for the Riesz transforms associated with the
Gaussian measure, {\em Rev. Mat. Iberoamericana} 10 (1994),
229--281.

\bibitem{FoSc} {L. Forzani and R. Scotto}, The higher order Riesz
transform for Gaussian measure need not be of weak type $(1,1)$,
{\em Studia Math.} 131 (1998), 205--214.

\bibitem{GaMaSjTo} {J. Garc\'{\i}a-Cuerva, G. Mauceri, P. Sj\"ogren and J.L
Torrea}, Higher-order Riesz operators for the Ornstein-Uhlenbeck
semigroup, {\em Potential Anal.} 10 (1999), 379--407.

\bibitem{GaStTr} {G. Gasper, K. Stempak and W. Trebels}, Fractional
integration for Laguerre expansions, {\em Methods Appl. Anal.} 2
(1995), 67--75.

\bibitem{GLLNU} {P. Graczyk, J.-L. Loeb, I.A. L\'opez, A. Nowak and W.
Urbina}, Higher order Riesz transforms, fractional derivatives and
Sobolev spaces for Laguerre expansions, {\em J. Math. Pures Appl.}
(9) 84 (2005), 375--405.

\bibitem{GuInTo} {C. Guti\'{e}rrez, A. Incognito and J.L. Torrea}, Riesz
tranforms, $g$-functions and multipliers for the Laguerre semigroup,
{\em Houston J. Math.} 27 (2001), 579--592.

\bibitem{GuSeTo} { C. Guti\'{e}rrez, C. Segovia and J.L. Torrea}, On higher Riesz
transforms for Gaussian measures, {\em J. Fourier Anal. Appl.} 2
(1996), 583--596.

\bibitem{HaRoSeTo} {E. Harboure, L. de Rosa, C. Segovia and J.L.
Torrea}, $L^p$- dimension free boundedness for Riesz transforms
associated to Hermite functions, {\em Math. Ann.} 328 (2004),
653--682.

\bibitem{HSTV} {E. Harboure, C. Segovia, J.L. Torrea and B.
Viviani}, Power weighted $L^p$-inequalities for Laguerre-Riesz
transforms, to appear in {\em Arkiv f\"or Matematik}.

\bibitem{HaToVi} {E. Harboure, J.L. Torrea and Viviani}, Riesz
transforms for Laguerre expansions, {\em Indiana Univ. Math. J.} 55
(2006), 999--1014.

\bibitem{Kanj} {Y. Kanjin and E. Sato}, The Hardy-Littlewood
theorem on fractional integration for Laguerre series, {\em Proc.
Amer. Math. Soc.} 123 (1995), 2165--2171.

\bibitem{Lebe} {N.N. Lebedev}, {\em Special functions and their
applications}, Dover Publications, Inc., New York, 1972.

\bibitem{Muck1} {B. Muckenhoupt}, Poisson integrals for Hermite
and Laguerre expansions, {\em Trans. Amer. Math. Soc.} 139 (1969),
231--242.

\bibitem{Muck2} {B. Muckenhoupt}, Hermite conjugate expansions,
{\em Trans. Amer. Math. Soc.} 139 (1969), 243--260.

\bibitem{Muck3} {B. Muckenhoupt}, Conjugate functions for Laguerre
expansions, {\em Trans. Amer. Math. Soc.} 147 (1970), 403--418.

\bibitem{Muck4} { B. Muckenhoupt}, Hardy's inequality with weights, {\em Studia Math.}
44 (1972), 31--38.

\bibitem{NoSt1} {A. Nowak and K. Stempak}, Weighted estimates for the Hankel transform
transplantation operator, {\em Tohoku Math. J.} 58 (2006),
277-301.

\bibitem{NoSt2} {A. Nowak and K. Stempak}, Riesz transforms and
conjugacy for Laguerre function expansions of Hermite type, {\em
J. Funct. Anal.} 244 (2007), 399-443.

\bibitem{PeSo} {S. P\'{e}rez and F. Soria}, Operators associated with the Ornstein-Uhlenbeck
semigroup, {\em J. London Math. Soc. (2)} 61 (2000), 857--871.

\bibitem{Sjog} {P. Sj\"ogren}, {\em On the maximal function for the Mehler
kernel}, Lecture Notes in Math. 992, Springer, Berlin, 1983.

\bibitem{Stei} {E. M. Stein}, {\em Topics in Harmonic
Analysis related to the Littlewood-Paley theory}, Ann. of Math.
Studies 63, Princenton Univ. Press, Princenton, 1970.

\bibitem{Stem} {K. Stempak}, Heat-diffusion and Poisson integrals
for Laguerre expansions, {\em Tohoku Math. J.(2)} 46 (1994),
83--104.

\bibitem{StTo1} {K. Stempak and J.L. Torrea}, Poisson integrals and
Riesz transforms for Hermite function expansions with weights, {\em
J. Functional Anal.} 202 (2003), 443--472.

\bibitem{StTo2} {K. Stempak and J.L. Torrea}, Higher Riesz transforms and
imaginary powers associated to the harmonic oscillator, {\em Acta
Math. Hungar.} 111 (2006), 43-64.

\bibitem{Szeg} {G. Szeg\"o}, {\em Orthogonal polynomials}, Colloquium Publications,
Vol. XXIII, American Math. Soc., Providence, R.I., 1975.

\bibitem{Than} {S. Thangavelu}, {\em Lectures on Hermite and Laguerre
expansions}, Mathematical Notes, 42, Princenton University Press,
Princenton, 1993.

\bibitem{Urbi} {W. Urbina}, On singular integrals with respect to
the Gaussian measure, {\em Ann. Scuola Norm. Sup. Pisa Cl. Sci. (4)}
17 (1990), 531--567.

\bibitem{Wats} {G. N. Watson}, {\em A treatise on the theory of
Bessel functions}, Cambridge Univ. Press, Cambridge, 1995.

\end{thebibliography}
\end{document}